\documentclass[a4paper,11pt]{article}
\usepackage{graphicx}
\usepackage{amsfonts}
\usepackage{amsmath}
\usepackage{amssymb}
\usepackage{indentfirst}
\usepackage{titlesec}
\usepackage{caption2}

\def\barr{\begin{array}}
\def\earr{\end{array}}
\def\bali{\begin{aligned}}
\def\eali{\end{aligned}}
\def\bearr{\begin{eqnarray}}
\def\eearr{\end{eqnarray}}

\providecommand{\play}{\displaystyle}

\providecommand{\li}{\limits}

\providecommand{\pt}{\partial}

\providecommand{\ra}{\rightarrow}

\providecommand{\da}{\downarrow}

\providecommand{\Prob}{\mathbf P}

\providecommand{\E}{\mathbf E}

\providecommand{\al}{\alpha}

\providecommand{\bt}{\beta}

\providecommand{\gm}{\gamma}

\providecommand{\Gm}{\Gamma}

\providecommand{\dt}{\delta}

\providecommand{\ve}{\varepsilon}

\providecommand{\tht}{\theta}

\providecommand{\kp}{\kappa}

\providecommand{\lb}{\lambda}

\providecommand{\Lb}{\Lambda}

\providecommand{\sm}{\sigma}

\providecommand{\Sm}{\Sigma}

\providecommand{\R}{\mathbb R}

\providecommand{\cE}{\mathcal E}

\providecommand{\cG}{\mathcal G}

\providecommand{\cS}{\mathcal S}

\providecommand{\iY}{\mathfrak{Y}}

\providecommand{\iiY}{\mathbf{\mathfrak{Y}}}

\providecommand{\grad}{\nabla}

\providecommand{\heff}{\mathbf{h}_{\text{eff}}}

\providecommand{\torus}{\mathbb{T}^2}

\providecommand{\graph}{\mathbb{G}}

\begin{document}

\author{Mark Freidlin\thanks{Dept of Mathematics, University of Maryland at
College Park, mif@math.umd.edu.}, Wenqing Hu\thanks{Dept of
Mathematics, University of Maryland at College Park,
huwenqing@math.umd.edu.}}
\title{On perturbations of generalized Landau-Lifshitz dynamics}
\date{}

\maketitle
\bigskip

\begin{abstract}

We consider deterministic and stochastic perturbations of dynamical
systems with conservation laws in $\R^3$. The Landau-Lifshitz
equation for the magnetization dynamics in ferromagnetics is a
special case of our system. The averaging principle is a natural
tool in such problems. But bifurcations in the set of invariant
measures lead to essential modification in classical averaging. The
limiting slow motion in this case, in general, is a stochastic
process even if pure deterministic perturbations of a deterministic
system are considered. The stochasticity is a result of
instabilities in the non-perturbed system as well as of existence of
ergodic sets of a positive measure. We effectively describe the
limiting slow motion.

\end{abstract}



\section{Introduction}

The analytical study of magnetization dynamics governed by the
Landau-Lifshitz equation (see [18]) has been the focus of
considerable research for many years. In normalized form this
equation reads as (see [5], equations (2.51) and (2.53)):

$$\dfrac{\pt\mathbf{m}}{\pt t}=-\mathbf{m}\times \heff  - \al \mathbf{m} \times (\mathbf{m} \times \heff)
\ , \ \mathbf{m}(\mathbf{r},0)=\mathbf{m}_0(\mathbf{r})\in \R^3 \ ,
\ |\mathbf{m}_0(\mathbf{r})|=1 \ . \eqno(1.1)$$

Here $\heff$ is an effective field. The three-dimensional vector
$\mathbf{m}(\mathbf{r},t)$ is the magnetization of the material at
 a fixed point $\mathbf{r}\in \R^3$ at time $t$; The term $\al
\mathbf{m} \times (\mathbf{m} \times \heff)$ is the Landau-Lifshitz
damping term, $0<\al<<1$. One can check that (1.1) preserves a first
integral $F(\mathbf{m})=\dfrac{1}{2}|\mathbf{m}|^2$. Therefore for
fixed $\mathbf{r}$, the system (1.1) describes a motion on the
sphere in $\R^3$.

One can introduce an energy density function $G$ such that $\grad G
= -\heff$. Then equation (1.1) can be written as follows:

$$\dfrac{d\mathbf{m}}{dt}=\mathbf{m}\times \grad G  + \al \mathbf{m} \times (\mathbf{m} \times \grad G)
\ , \ \mathbf{m}(0)=\mathbf{m}_0\in \R^3 \ , \ |\mathbf{m}_0|=1 \ .
\eqno(1.2)$$

We assume that $G$ is a smooth generic function. Considered on the
unit sphere $S^2$ in $\R^3$, such a function may have three types of
critical points: maxima, minima and saddle points. Without the
damping $\al \mathbf{m} \times (\mathbf{m} \times \grad G)$ the
energy density $G$ is preserved. One easily checks that $\grad
G\cdot (\al \mathbf{m} \times (\mathbf{m} \times \grad
G))=-\al|\mathbf{m}\times \grad G|^2$ so that the damping term is a
kind of "friction" for the system (1.2), just like the classical
friction in Hamiltonian systems (compare with [4]).

If $0<\al<<1$, the dynamics of (1.2) has two distinct time scales:
the fast time scale of the precessional dynamics and the relatively
slow time scale of relaxational dynamics caused by the small damping
term $\al \mathbf{m}\times (\mathbf{m} \times \grad G)$. Therefore
it is natural to use the averaging principle to describe the
long-time evolution of energy density $G$. However the classical
averaging principle here should be modified: existence of saddle
points of $G(\mathbf{m})$ on the sphere $\{|\mathbf{m}|=1\}$ leads
to stochastic, in a certain sense, behavior of the slow motion even
in the case of purely deterministic damping term (compare with [4]).
Moreover, in Section 5, we consider a more general class of
equations, where level set components of first integrals, which are
compact two-dimensional surfaces may have topological structure
different from a sphere. If genus of such a surface is positive, the
non-perturbed system can have positive area ergodic sets. Existence
of such sets lead to an "additional stochasticity". Description of
the stochastic process which characterizes the long-time evolution
of the energy is one of the main goals of this paper.

Random perturbation caused by thermal fluctuations become
increasingly pronounced in nano-scale devices. To take this into
account one can include in the right-hand side of (1.2) a small
stochastic term. This stochastic term, in general, introduces one
more time scale in the system. Interplay between the influence of
small damping and even smaller stochastic term leads to certain
changes in the metastability of the system. Description of the
metastable distributions is another goal of this paper. There are
some other asymptotic regimes of the Landau-Lifshitz dynamics which
we mention briefly and we will consider them in more details
elsewhere.

\

\section{Sketch of the paper}

In this section we give an informal sketch of the results.

In the next two sections we consider perturbations of the following
equation

$$\dot{\widetilde{X}}_t=\grad F (\widetilde{X}_t)\times
\grad G (\widetilde{X}_t) \ , \ \widetilde{X}_0=x_0\in \R^3 \ ,
\eqno(2.1)$$ which could be regarded as a generalized
Landau-Lifshitz equation.

Here $G(x)$ and $F(x)$, $x\in \R^3$, are smooth enough generic
functions (this means that each of these functions has a finite
number of critical points which are assumed to be non-degenerate),
$\lim\li_{|x| \ra \infty}F(x)=\infty$. The initial point
$x_0=x_0(z)$ is chosen in such a way that $F(x_0(z))=z$. As before
we call $G(x)$ energy (to be precise, $G(x)$ in (1.2) is the energy
density but for brevity we call it energy).

It is easy to see that $F(x)$ and $G(x)$ are first integrals of
system (2.1). For instance,

$$\dfrac{dF(\widetilde{X}_t)}{dt}=\grad F (\widetilde{X}_t)\cdot
(\grad F (\widetilde{X}_t)\times \grad G (\widetilde{X}_t))=0 \ .$$

Note also that the Lebesgue measure in $\R^3$ (the volume) is
invariant for system (2.1): $$\text{div} (\grad F (x) \times \grad G
(x))=\grad G (x) \cdot (\grad \times \grad F (x))-\grad F(x) \cdot
(\grad \times \grad G(x))=0 \ .$$ This implies, in particular, that
$\play{\dfrac{1}{|\grad F(x)|}}$ is the density of an invariant
measure of system (2.1) considered on the surface
$\widetilde{S}(z)=\{x\in \R^3: F(x)=z\}$ with respect to the area on
$\widetilde{S}(z)$. Notice that the surface $\widetilde{S}(z)$ may
have several connected components. For brevity in the next two
sections, and in the rest of this section (except the last four
paragraph), we assume that the level surface
$\widetilde{S}(z)=\{x\in \R^3: F(x)=z\}$ has only one connected
component and this component is homeomorphic to $S^2$. In Sections 5
and 6 we will drop this assumption and consider more general
situations.

As we already mentioned, the damping term in (1.2) preserves the
first integral $\play{\dfrac{1}{2}|\mathbf{m}|^2}$, so that we
consider, first, perturbations of (2.1) preserving $F(x)$. The
perturbed equation can be written in the form

$$\dot{\widetilde{X}}_t^{\ve}=\grad F (\widetilde{X}_t^{\ve})\times \grad G (\widetilde{X}_t^\ve)
+\ve \grad F (\widetilde{X}_t^\ve)\times
\widetilde{\mathbf{b}}(\widetilde{X}_t^\ve) \ , \
\widetilde{X}_0^\ve=x_0\in \R^3 \ . \eqno(2.2)$$

Here $\widetilde{\mathbf{b}}(\bullet)$ is a smooth vector field in
$\R^3$. In the next two sections we assume for brevity that the
perturbation $\ve \grad F \times \widetilde{\mathbf{b}}$ is of
"friction" type:
$$\grad G(x) \cdot (\grad F(x) \times \widetilde{\mathbf{b}}(x))<0 \
, \ x \in \widetilde{S}(z)\subset \R^3. \eqno(2.3)$$

Note that any vector field $\grad F(x) \times
\widetilde{\mathbf{b}}(x)$ can be written in the form $\grad F(x)
\times (\grad F(x) \times \mathbf{b}(x))$ for some vector field
$\mathbf{b}(x)\in \R^3$. Indeed, without loss of generality one can
assume that $\widetilde{\mathbf{b}}(x)\perp \grad F$. Each such
vector $\widetilde{\mathbf{b}}(x)$ can be represented as $\grad F
(x) \times \mathbf{b}(x)$. So that the perturbed equation can be
written as

$$\dot{\widetilde{X}}_t^\ve=\grad F(\widetilde{X}_t^\ve) \times \grad G(\widetilde{X}_t^\ve)+
\ve \grad F(\widetilde{X}_t^\ve)\times (\grad
F(\widetilde{X}_t^\ve)\times \mathbf{b}(\widetilde{X}_t^\ve)) \ ,
\widetilde{X}_0^\ve=x_0\in \R^3 \ . \eqno(2.4)$$

Furthermore, using the identity $\mathbf{A}\cdot(\mathbf{B}\times
(\mathbf{C} \times \mathbf{D}))=(\mathbf{A}\times \mathbf{B})\cdot
(\mathbf{C}\times \mathbf{D})$ we can check that

$$\grad G \cdot (\grad F \times (\grad F \times \mathbf{b}))=-(\grad F \times
\mathbf{b})\cdot (\grad F \times \grad G) \ . \eqno(2.5)$$

Therefore the "friction-like" condition (2.3) becomes

$$(\grad F \times \mathbf{b})\cdot (\grad F \times \grad G)>0  \ . \eqno(2.6)$$

The equation (1.2) corresponds to the case that
$\mathbf{b}(\widetilde{X}_t^\ve)=\grad G (\widetilde{X}_t^\ve)$ and
$F(\widetilde{X}_t^\ve)=\dfrac{1}{2}|\widetilde{X}_t^\ve|^2$. One
easily checks that system (2.4) preserves $F$ so that
$\widetilde{X}_t^\ve$ is moving on a certain level surface
$\{F=z\}$.

We make some geometric assumptions that are used in Sections 3 and
4. Suppose that the set $S(z)=\{x\in \R^3: G(x)\leq
G(x_0(z))+1\}\cap \{x\in \R^3: F(x)=z\}$ is a 2-dimensional
Riemannian manifold which is $C^{\infty}$-diffeomorphic to
$R=\{(a,b)\in \R^2: a^2+b^2\leq 1\}$. Let the $C^\infty$
diffeomorphism be $f: S(z)\ra R$. To be specific, we denote
$f(x_1,x_2,x_3)=(f_1(x_1,x_2,x_3),f_2(x_1,x_2,x_3))$ for
$(x_1,x_2,x_3)\in S(z)$. We assume that the diffeomorphism $f$ is
non-singular for $(x_1,x_2,x_3)\in S(z)$. We denote by $d(\bullet,
\bullet)$ the metric on $S(z)$ induced by standard Euclidean metric
in $\R^3$. Let our function $G$ on $S(z)$ have only one saddle point
and two minima, and these critical points are non-degenerate. Assume
that the level surfaces $\{G=g\}$ are transversal to the level
surface $\{F=z\}$: $\grad F(x)$ and $\grad G (x)$ are not parallel.
We denote by $C(g,z)$ the set $\{G(x)=g\}\cap \{F(x)=z\}$. Without
loss of generality we can assume that $C(0,z)=\{G=0\}\cap \{F=z\}$
is the $\infty$-shaped curve (homoclinic trajectory) on $\{F=z\}$
corresponding to the saddle point of $G$. Let the saddle point of
$G$ on $\{F=z\}$ be $O_2(z)$ and the two minima be $O_1(z)$ and
$O_3(z)$. Suppose that as $z$ varies, the curves $O_1(z), O_2(z)$
and $O_3(z)$ are transversal to $\{F=z\}$ (see Fig.1). Notice that
when $g>0$, $C(g,z)$ has only one connected component which we call
$C_2(g,z)$. When $g<0$, $C(g,z)$ has two connected components
$C_1(g,z)$ and $C_3(g,z)$ bounding domains on $S(z)$ containing
$O_1(z)$ and $O_3(z)$ respectively. Let $C_1(0,z)$ and $C_3(0,z)$ be
the parts of the homoclinic trajectory $C(0,z)$ bounding domains
containing $O_1(z)$ and $O_3(z)$ respectively. Let
$C_2(0,z)=C(0,z)$. Let $D_i(g,z)$ $(i=1,2,3)$ be the region bounded
by $C_i(g,z)$.

In the next two sections when we speak about a stochastic process or
a motion on the surface $S(z)$, for example $X_t^\ve$,
$X_t^{\ve,\dt}$ etc. , we are assuming that they are stopped once
they hit $\partial S(z)$.

\begin{figure}
\centering
\includegraphics[height=6.5cm, width=10cm , bb=81 13 433 235]{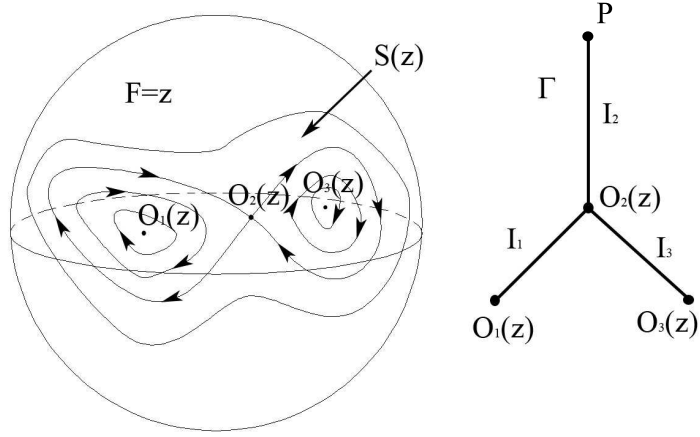}
\caption{The Landau-Lifshitz dynamics}
\end{figure}

To study equation (2.3), we make a time change $t\mapsto
\dfrac{t}{\ve}$. Let $X_t^\ve=\widetilde{X}_{t/\ve}^\ve$. We get
from (2.2) that

$$\dot{X}_t^\ve=\dfrac{1}{\ve}\grad F(X_t^\ve) \times \grad G(X_t^\ve)+
\grad F(X_t^\ve)\times (\grad F(X_t^\ve)\times \mathbf{b}(X_t^\ve))
\ , X_0^\ve=x_0\in \R^3 \ , 0<\ve<<1 \ . \eqno(2.7)$$

Therefore the fast motion is defined by the vector field
$\dfrac{1}{\ve}\grad F \times \grad G$, and the slow motion is due
to $\grad F \times (\grad F \times \mathbf{b})$ (we will sometimes
ignore the arguments since they could be directly understood from
the context). In order to study the limiting behavior of the process
$X_t^\ve$, we introduce a graph $\Gm$ (compare with [14, Chapter
8]). The graph $\Gm$ is constructed in the following way. Let us
identify the points of each connected component of the level sets of
$G$ on $S(z)$. Let the identification mapping be $\iY$. The set
obtained after such an identification, equipped with the natural
topology, is a graph $\Gm$ with an interior vertex $O_2(z)$
corresponding to the saddle point $O_2(z)$ on $S(z)$ and related
homoclinic curve (in the following we will use the same symbol for
either the critical point of $G$ on $S(z)$ or the corresponding
vertex on $\Gm$), and two exterior vertices $O_1(z)$ and $O_3(z)$
corresponding to the stable equilibriums $O_1(z)$ and $O_2(z)$ on
$S(z)$, together with another exterior vertex $P$ corresponding to
$\pt S(z)$ (notice that by our definition $S(z)=\{x\in \R^3:
G(x)\leq G(x_0(z))+1\}\cap \{x\in \R^3: F(x)=z\}$ so that $\pt S(z)$
is a level curve of $G$ on $\{F=z\}$). The edges of the graph are
defined as follows: edge $I_2$ corresponds to trajectories on $S(z)$
lying outside $C(0,z)$; edges $I_1$ and $I_3$ correspond to those
trajectories on $S(z)$ belonging to the wells containing $O_1(z)$
and $O_3(z)$, respectively. A point $\iY(x)=y\in \Gm$ can be
characterized by two coordinates $(g,k)$ where $g=G(x)$ is the value
of function $G$ at $x\in \iY^{-1}(y)\subset S(z)$, and $k=k(x)$ is
the number of the edge of the graph $\Gm$ to which $y=\iY(x)$
belongs. Notice that $k$ is not chosen in a unique way since for
$y=O_2(z)$ the value of $k$ can be either $1$, $2$ or $3$. The
distance $\rho(y_1,y_2)$ between two points $y_1=(G(x_1),k)$ and
$y_2=(G(x_2),k)$ is simply $\rho(y_1,y_2)=|G(x_1)-G(x_2)|$. For
$y_1, y_2\in \Gm$ belonging to different edges of the graph it is
defined as $\rho(y_1,y_2)=\rho(y_1,O_2(z))+\rho(O_2(z),y_2)$.

The slow component of $X_t^\ve$ is the projection of $X_t^\ve$ on
$\Gm$: $Y_t^\ve=\iY(X_t^\ve)$. Using the classical averaging
principle one can describe the limiting motion of $Y_t^\ve$ as $\ve
\da 0$ inside the edges. But it turns out that the trajectory
$Y_t^\ve$, when hitting the interior vertex $O_2(z)$ on $\Gm$, is
very sensitive to $\ve$. This means that $Y_t^\ve=\iY(X_t^\ve)$,
$G(X_0^\ve)>G(O_2(z))$, hits $O_2(z)$ in a finite time $t_0^\ve$
such that $\lim\li_{\ve \da 0}t_0^\ve=t_0$ exists and finite, and
after that alternatively as $\ve \da 0$ goes to $I_1$ or $I_3$. The
limit of $Y_t^\ve$ as $\ve \da 0$ for $t>t_0$ does not exist
(compare with [4]). In order to describe the limiting behavior, we
have to regularize the problem. To this end one can add a small
stochastic perturbation of order $\dt$ either to the initial
condition or to the equation. Let $X_t^{\ve,\dt}$ be the result of
addition of such a perturbation. Then, under certain mild
assumptions, the slow component $\iY(X_t^{\ve,\dt})$ of
$X_t^{\ve,\dt}$ converges weakly in the space of continuous
trajectories on any finite time interval $[0,T]$ to a stochastic
process $Y_t$ on the graph $\Gm$ as first $\ve \da 0$ and then $\dt
\da 0$. Since small random perturbations, as a rule, are available
in the system, exactly this weak limit characterizes the behavior of
$\widetilde{X}_{t/\ve}^\ve$ as $0<\ve<<1$. We will introduce
different types of regularization and prove that all these
regularizations lead to the same limiting stochastic process $Y_t$
on $\Gm$, which we calculate.

The proofs, in Section 3 and, partly, in Section 4, are similar to
the case of perturbations of Hamiltonian systems ([14, Chapter 8],
[4]), and we pay most of the attention to the arguments which are
not presented in these works. For instance, in the case of
regularization by a random perturbation of the initial point, bounds
for the hitting time of the homoclinic trajectory are considered in
details.

So far we considered just deterministic perturbations preserving the
first integral $F$. Stochastic perturbations were used just for
regularization of the problem. One can consider also
white-noise-type perturbations preserving $F$ of the same or of a
larger order than deterministic perturbations. Then, in an
appropriate time scale, the limiting slow motion converges to a
diffusion process on a graph (Section 4). In general, deterministic
and stochastic perturbations have different order, so that, after
time rescaling $t \ra \dfrac{t}{\ve}$, the perturbed equation has
the form

$$X_t^{\ve,\dt}=\dfrac{1}{\ve}\grad F (X_t^{\ve,\dt})\times \grad G (X_t^{\ve,\dt})+\grad F (X_t^{\ve,\dt})
\times \widetilde{\mathbf{b}}(X_t^{\ve,\dt})+
\dt\sm(X_t^{\ve,\dt})\circ \dot{W}_t\ , X_0^{\ve,\dt}=x_0 \ .
\eqno(2.8)$$

Here $\dot{W}_t$ is the standard Gaussian white noise, $\sm(x)$ is a
smooth matrix-function such that $\sm^T\grad F \equiv 0$. If we
denote by $a(x)=\sm(x)\sm^T(x)$ the diffusion matrix, the condition
$\sm^T(x)\grad F(x)\equiv 0$ is equivalent to the assumption that
$a(x)\grad F(x)\equiv 0$. The stochastic term in (2.8) is understood
in the Stratonovich sense, then $F(X_t^{\ve,\dt})\equiv F(x_0)$ with
probability 1. We assume that the matrix $a$ is non-degenerate on
$\{F=z\}$. (We will specify the non-degeneracy in Section 4.)

The process $X_t^{\ve,\dt}$ defined by (2.8) lives on the surface
$\{x\in \R^3: F(x)=z\}$ and has a slow and a fast component as
$\ve<<1$ and $\dt
> 0 $ fixed. The slow component is again the projection
$\iY(X_t^{\ve,\dt})$ of $X_t^{\ve,\dt}$ on the graph $\Gm$. We
consider the case $0<\ve<<\dt<<1$ and assume that the deterministic
perturbation is friction-like.

If $\ve>0$ is small enough and $\dt=0$, the system $X_t^{\ve,0}$,
$X_0^{\ve,0}=x\in S(z)$, has three critical points $O_1'(z)$,
$O_2'(z)$, $O_3'(z)$ of the same type as the corresponding points
$O_i(z)$. The distance between corresponding points tends to zero
together with $\ve$. If $0<\dt<<1$, $X_t^{\ve,\dt}$,
$X_0^{\ve,\dt}=x$, at a time $t=T^\dt(\lb)$, $\lim\li_{\dt \da
0}\dt^2 \ln T^\dt(\lb)=\lb>0$, is situated in a small neighborhood
of the metastable state $M^\ve(x,\lb)$; $M^\ve(x,\lb)$ is one of the
stable equilibriums of $X_t^{\ve,0}$. The function $M^\ve(x,\lb)$ is
defined by the action functional for the family $X_t^{\ve,\dt}$ as
$\dt \da 0$ (see [8], [10], [13], [14], [21]).

But if $\ve$ tends to zero, the situation is different:
$X_{T^\dt(\lb)}^{\ve,\dt}$, $X_0^{\ve,\dt}=x_0$, converges to a
random variable distributed between $O_1(z)$ and $O_3(z)$ as
$0<\ve<<\dt<<1$. The set of possible distributions between the
minima is finite and is independent of the stochastic part of
perturbations. But which of these distributions is realized at a
time $T^\dt(\lb)$ depends on $\lb$ and $x_0=X_0^{\ve,\dt}$, as well
as on stochastic perturbations. We describe these metastable
distributions in Section 4.

\

Perturbations of a more general equation than (2.1) are considered
in Section 5. The non-perturbed motion in this case, in general, has
just one smooth first integral and the averaging procedure
essentially depends on the topological structure of the connected
components of level sets of the existing first integral. Each
connected component is two dimensional orientable compact manifold.
The topology of such a manifold is determined by its genus. We show
that if the genus is greater than zero (for instance, when the
component is a 2-torus $\torus$), the limiting slow motion spends an
exponentially distributed random time at some vertices.

Perturbations of system (2.1) may have different origin and they may
have different order. In the last Section 6, we briefly consider
such a situation.

Perturbations of (2.1) breaking both first integrals $F(x)$ and
$G(x)$ can be considered: (after time change)

$$\dot{X}_t^\ve=\dfrac{1}{\ve}\grad F(X_t^\ve) \times \grad G(X_t^\ve)+
 \mathbf{B} (X_t^\ve) \ , X_0^\ve=x_0(z)\in \R^3 \ ,
0<\ve<<1 \ . \eqno(2.9)$$ Here $\mathbf{B}(\bullet)$ is a general
smooth vector field on $\R^3$. Then the perturbed motion is not
restricted to the level surface $\{F=z\}$. In this case the slow
component of the perturbed motion lives on an "open book" $\sqcap$
homeomorphic to the set of connected components of the level sets
$C(z_1,z_2)=\{x\in \R^3: F(x)=z_1, G(x)=z_2\} \ , \ (z_1, z_2) \in
\R^2$ (compare with [16]). The slow component of the motion is equal
to $\iiY(X_t^\ve)=Y_t^\ve$, where $\iiY: \R^3 \ra \sqcap$ is the
identification mapping. After an appropriate regularization,
$Y_t^\ve$ approaches as $\ve \da 0$ a stochastic process $Y_t$ on
$\sqcap$. We will consider this question in more details elsewhere.

\section{Regularization by perturbation of the initial condition}

We study in this section the regularization of system (2.7) by a
stochastic perturbation of the initial condition.

Let $U_\dt(x)=\{y\in S(z): d(x,y)<\dt\}$.

Consider the equation:

$$\dot{X}_t^{\ve,\dt}=\dfrac{1}{\ve}\grad F(X_t^{\ve,\dt}) \times \grad G(X_t^{\ve\,\dt})+
\grad F(X_t^{\ve,\dt})\times (\grad F(X_t^{\ve,\dt})\times
\mathbf{b}(X_t^{\ve,\dt})) \ , X_0^{\ve,\dt}=x_0(z,\dt)\in \R^3 \ .
\eqno(3.1)$$

Here $0<\dt<<1$ is a small parameter. The initial position
$x_0(z,\dt)=X_0^{\ve,\dt}$ is a random variable distributed
uniformly in $U_{\dt}(x_0(z))\subset \{F=z\}$. We are choosing $\dt$
small enough so that $U_\dt(x_0(z))\subset S(z)$.

Our goal is to prove the following

\

\textbf{Theorem 3.1.}  \textit{Let $X_t^{\ve,\dt}$ be the solution
of equation (3.1), and $Y_t^{\ve,\dt}=\iY(X_t^{\ve,\dt})$ be the
slow component of $X_t^{\ve,\dt}$. Then, for each $T>0$,
$Y_t^{\ve,\dt}$ converges weakly in the space of continuous
functions $f$: $[0,T]\ra \Gm$ to a stochastic process
$\overline{Y}_t(x_0(z))$ as, first, $\ve \da 0$ and then $\dt \da
0$.}

\

We will define the process $\overline{Y}_t(x_0(z))$ later in this
section.

\

Let us start with the perturbed, but not regularized system (2.7).
The motion of $X_t^\ve$ is on the surface $S(z)$. The change of
$G(X_t^\ve)$ is governed by the equation

$$\begin{array}{ll}
\dfrac{dG(X_t^\ve)}{dt}&=\grad G \cdot (\dfrac{1}{\ve}\grad F \times
\grad G+ \grad F \times (\grad F \times \mathbf{b}))\\&= \grad G
\cdot (\grad F \times (\grad F \times \mathbf{b}))\\&=-(\grad F
\times \mathbf{b}) \cdot (\grad F \times \grad G)  \ .
\end{array} $$

The function $G$ is a first integral of the unperturbed system (2.1)
and the damping term $\grad F \times (\grad F \times \mathbf{b})$ of
(2.7) plays the role of "friction" which makes the value of $G$
smaller and smaller.

The stable, but not asymptotically stable equilibriums $O_1(z)$ and
$O_3(z)$ of (2.1) become asymptotically stable equilibriums
$O_1'(z)$ and $O_3'(z)$ for the perturbed system (2.7). The saddle
point $O_2(z)$ becomes the saddle point $O_2'(z)$. The distances
between $O_1(z)$ ($O_2(z), O_3(z)$) and $O_1'(z)$ ($O_2'(z),
O_3'(z)$) are less than $A\ve$ for a constant $A>0$. When $\ve$ is
small enough, the pieces of the curves formed by $O_1'(z)$,
$O_2'(z)$ and $O_3'(z)$ (as $z$ varies) are transversal to
$\{F=z\}$. Separatrices of the saddle point $O_2'(z)$ are shown in
Fig.2. They, roughly speaking, divide the part of the surface $S(z)$
outside the $\infty$-shaped curve $C(0,z)$ in ribbons: the gray
ribbon enters the neighborhood of $O_1'(z)$, and the white ribbon
enters the neighborhood of $O_3'(z)$. The width of each ribbon is of
order $\ve$ as $\ve \da 0$.

\begin{figure}
\centering
\includegraphics[height=8cm, width=12cm , bb=62 50 381 257]{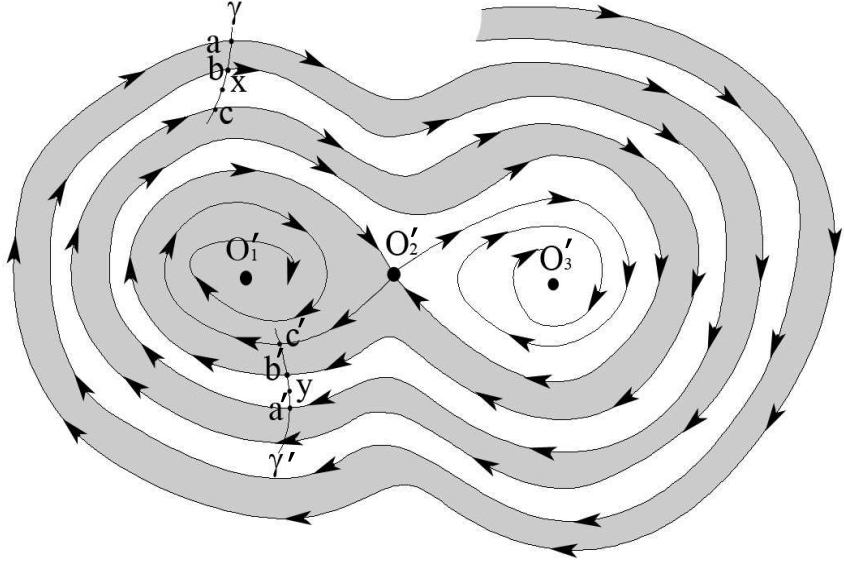}
\caption{White and grey ribbons}
\end{figure}

The trajectory $X_t^\ve$ has a fast component, which is close to the
non-perturbed motion (2.1) (with the speed of order
$\dfrac{1}{\ve}$), and the slow component, which is the projection
$Y_t^\ve=\iY(X_t^\ve)$ of $X_t^\ve$ on the graph $\Gm$ corresponding
to $G(x)$. Within each edge of the graph, say edge $I_i$, $i=1,2,3$,
standard averaging principle works. Let $G^\ve_t=G(X_t^\ve)$. We
have, by the standard averaging principle (cf. [1], Ch.10),
$$\lim\li_{\ve\da 0}\sup\li_{0\leq t \leq T<\infty}|G_t^\ve-G_t|=0 \
.$$ The function $G_t$ satisfies $G_0=G(x_0(z))$ and
$$\begin{array}{l}
\play{\dfrac{dG_t}{dt}}\play{=B^{(i)}(G_t)} \ ,  \text{ where } \\
\\
\play{B^{(i)}(g)}\\
\play{=\dfrac{1}{T_i(g)}\oint_{C_i(g,z)} \grad G \cdot (\grad F
\times (\grad F \times \mathbf{b}))\dfrac{dl}{|\grad F \times
\grad G|}}\\
\play{=-\dfrac{1}{T_i(g)}\oint_{C_i(g,z)}(\grad F \times
\mathbf{b})\cdot\dfrac{\grad F \times \grad G}{|\grad F \times \grad
G|}dl}\\
\play{=-\dfrac{1}{T_i(g)}\oint_{C_i(g,z)}(\grad F \times
\mathbf{b})\cdot
\mathbf{v} dl}\\
\play{=-\dfrac{1}{T_i(g)}\iint_{D_i(g,z)}\grad \times (\grad F
\times \mathbf{b})\cdot \mathbf{n} dm \ .}
\end{array} \eqno(3.2)$$

Here $\play{T_i(g)=\oint_{C_i(g,z)} \dfrac{dl}{|\grad F \times \grad
G|}}$ is the period of rotation for the unperturbed system (2.1)
along the curve $C_i(g,z)$. The vector $\textbf{v}=\dfrac{\grad
F\times \grad G}{|\grad F \times \grad G|}$ is the unit velocity
vector for the unperturbed system (2.1);
$\textbf{n}=\textbf{n}(x)=\dfrac{\grad F(x)}{|\grad F(x)|}$ is
normal to the level surface $\{F=z\}$. The area element on $\{F=z\}$
is denoted by $dm$. We used the Stokes formula in the last step.

Fix a point $x_0(z)$ on the level surface $\{F=z\}$ outside the
$\infty$ - shaped curve $C(0,z)$. To be specific, let $x_0(z)$
belong to the white ribbon. Let $\gm_s(z)$ be the curve on $\{F=z\}$
containing $x_0(z)$ and orthogonal to the perturbed trajectories
(2.7). Let $a(z),b(z),c(z)$ be the intersection points of $\gm_s(z)$
with separatrices neighboring to $x_0(z)$. To be specific, let
$x_0(z)$ lie between $b(z)$ and $c(z)$ (see Fig.3, where a part of
the flow is shown). By our transversality condition, we can take
$\lb>0$ small enough and a curve $\xi(\widetilde{z}) \ , \
\widetilde{z}\in [z-2\lb, z+2\lb]$ which lies on the surface
$\{G=G(x_0(z))\}$ and is transversal to the level surface $\{F=z\}$,
containing the point $x_0(z)$ ($\xi(z)=x_0(z)$). Let
$x_0(\widetilde{z})=\xi(\widetilde{z})$. Consider the curve
$\gm_s(\widetilde{z})$ on $\{F=\widetilde{z}\}$ containing
$x_0(\widetilde{z})$ and orthogonal to the trajectories of (2.7). We
also consider corresponding neighboring points $a(\widetilde{z})$,
$b(\widetilde{z})$, $c(\widetilde{z})$ defined for
$x_0(\widetilde{z})$ in the same way as we did for $x_0(z)$. For
fixed $\ve>0$, we choose $\lb$ small enough such that as
$\widetilde{z}$ varies in $[z-2\lb, z+2\lb]$, the curves
$a(\widetilde{z})$, $b(\widetilde{z})$ and $c(\widetilde{z})$ are
transversal to $\{F=z\}$. The part of $\gm_s(\widetilde{z})$ between
$a(\widetilde{z})$ ($b(\widetilde{z})$) and $b(\widetilde{z})$
($c(\widetilde{z})$) belongs to the grey (white) ribbon for the
trajectories of (2.7) on $\{F=\widetilde{z}\}$. Now we consider the
curvilinear rectangle $\square_1$ with vertices $a(z+\lb) \ ,\
a(z-\lb)\ , \ b(z-\lb) \ , \ b(z+\lb)$ constructed in the following
way: $\square_1$ consists of the parts of the curves of
$\gm_s(\widetilde{z})$ from $a(\widetilde{z})$ to $b(\widetilde{z})$
as $\widetilde{z}$ varies in $[z-\lb, z+\lb]$. We construct another
curvilinear rectangle $\square_2$ with vertices $b(z+\lb)\ , \
b(z-\lb) \ , \ c(z-\lb) \ , \ c(z+\lb)$ in exactly the same way as
$\square_1$, but consisting of curves $\gm_s(\widetilde{z})$ from
$b(\widetilde{z})$ to $c(\widetilde{z})$ as $\widetilde{z}$ varies
in $[z-\lb, z+\lb]$.

\begin{figure}
\centering
\includegraphics[height=9cm, width=12cm , bb=43 9 412 290]{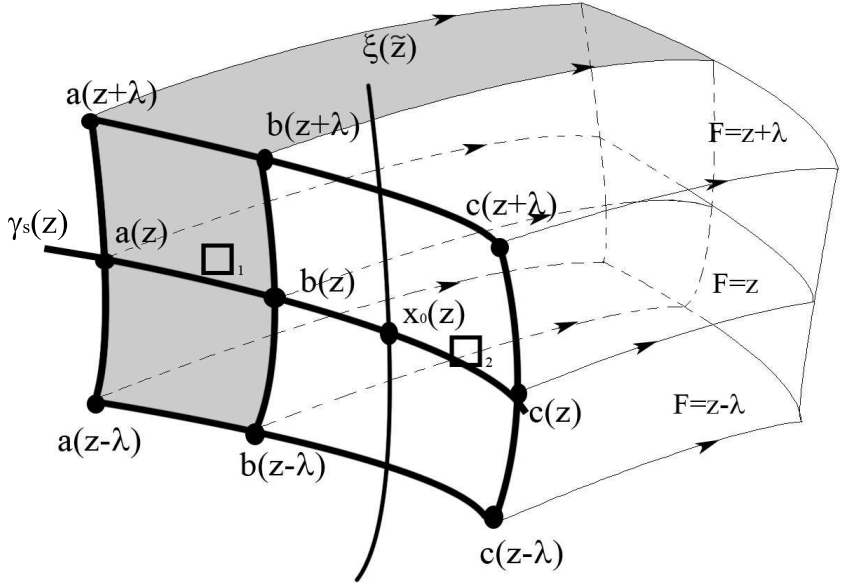}
\caption{Transversality}
\end{figure}

Let vector $\vec{\nu}$ be the unit vector outward normal to these
two curvilinear rectangles $\square_1$ and $\square_2$, pointing in
the direction opposite to the perturbed flow (2.7). By the
divergency theorem in $\R^3$ we see that, for $k=1,2$,

$$\begin{array}{l}
\play{\iint_{\square_k}\left(\dfrac{1}{\ve}\grad F \times \grad
G+\grad F \times (\grad F \times \mathbf{b})\right)\cdot
\vec{\nu}dm}
\\
\play{=\iiint_{\cE_k} \text{div}(\grad F \times (\grad F \times
\mathbf{b}))dV}
\\
\play{=-\int_{z-\lb}^{z+\lb}d\widetilde{z}\iint_{\cS_k(\widetilde{z})}\grad\times(\grad
F\times \mathbf{b})\cdot \dfrac{\grad F}{|\grad F|}dm}
\\
\play{=-2\lb \iint_{\cS_k(z)}\grad\times(\grad F\times
\mathbf{b})\cdot \dfrac{\grad F}{|\grad F|}dm+ o(\lb) \ .}
\end{array} \eqno(3.3)$$

We have used here the formula $\text{div}(\textbf{A}\times
\textbf{B})=\textbf{B}\cdot ( \grad \times \textbf{A} ) - \textbf{A}
\cdot (\grad \times \textbf{B})$. The regions $\cE_k$ and
$\cS_k(\widetilde{z})$ ($k=1,2$) are defined as follows: $\cE_k$ is
the 3-dimensional region filled by trajectories of (2.7) starting
from $\square_k$ and belonging to the family of level surfaces
$\{F=\widetilde{z}\}$, $\widetilde{z}\in [z-\lb, z+\lb]$, $k=1,2$;
$\cS_k(\widetilde{z})$ is the 2-dimensional region filled by
trajectories of (2.7) starting from
$\square_k\cap\{F=\widetilde{z}\}$ and restricted to the family of
level surfaces $\{F=\widetilde{z}\}$, $k=1,2$, $\widetilde{z}\in
[z-\lb, z+\lb]$. Notice that the boundary of the compact set $\cE_k$
consist of $\square_k$ and a surface formed by the perturbed
trajectory. For notational convenience the area element on
$\square_k$ ($k=1,2$) is denoted also by $dm$.

Let $L(a(z),b(z))$ and $L(b(z),c(z))$ be, respectively, the arc
length of $\gm_s(z)$ between $a(z)$ and $b(z)$, and between $b(z)$
and $c(z)$. The flux of the vector field $\play{\dfrac{1}{\ve}\grad
F \times \grad G + \grad F \times (\grad F \times \mathbf{b})}$
through $\square_k$ ($k=1,2$) is equal to
$\play{-\iint_{\square_k}\left|\dfrac{1}{\ve}\grad F \times \grad G
+\grad F \times (\grad F\times \mathbf{b})\right|dm}$.

Let $\text{Area}(\bullet)$ denote the area of some domain. Let
$|J(\widetilde{z},\gm_s(\widetilde{z}))|\neq 0$ be the Jacobian
factor between the area element on $\square_1\cup \square_2$ and
$d\widetilde{z}d\gm_s(\widetilde{z})$. We have

$$\begin{array}{l}
\text{Area}(\square_1\cup \square_2)\\
\play{=\int_{z-\lb}^{z+\lb}d\widetilde{z}\int_{a(\widetilde{z})}^{c(\widetilde{z})}
|J(\widetilde{z},\gm_s(\widetilde{z}))|d\gm_s(\widetilde{z})} \\
\play{=2\lb \int_{a(z)}^{c(z)} |J(z,\gm_s(z))|d\gm_s(z)+(I)}\\
\play{=2\lb |J(z,b(z))| L(a(z),c(z))+2\lb(II)+(I)}\ .
\end{array}$$

Here

$$(I)=\int_{z-\lb}^{z+\lb}d\widetilde{z}\left(\int_{a(\widetilde{z})}^{c(\widetilde{z})}
|J(\widetilde{z},\gm_s(\widetilde{z}))|d\gm_s(\widetilde{z})-\int_{a(z)}^{c(z)}
|J(z,\gm_s(z))|d\gm_s(z)\right) \ ,$$ and
$$(II)=\int_{a(z)}^{c(z)}(J(z,\gm_s(z))-J(z,b(z)))d\gm_s(z) \ .$$

Note that $|(I)|\leq C_1 \lb^2$ since the function
$I(\widetilde{z})=\play{\int_{a(\widetilde{z})}^{c(\widetilde{z})}
|J(\widetilde{z},\gm_s(\widetilde{z}))|d\gm_s(\widetilde{z})}$
satisfies $|I(\widetilde{z}_1)-I(\widetilde{z}_2)|\leq C_2
|\widetilde{z}_1-\widetilde{z}_2|$. We also have $|(II)|\leq
C_3L(a(z),c(z))^2$ since $|J(z,\gm_s(z))-J(z,b(z))|\leq C_4
|\gm_s(z)-b(z)|\leq C_5 L(a(z),c(z))$. Combining these estimates
with (3.3) and the fact that for some constants $C_6, C_7>0$,
$$\dfrac{C_6}{\ve}\leq\dfrac{1}{\text{Area}(\square_1\cup \square_2)}
\iint_{\square_1\cup \square_2}\left|\dfrac{1}{\ve}\grad F \times
\grad G + \grad F \times (\grad F \times
\mathbf{b})\right|dm\leq\dfrac{C_7}{\ve} \ , $$ we see that as
$\ve\da 0$, the asymptotic widths of the grey and white ribbons
(i.e. $L(a(z),b(z))$ and $L(b(z),c(z))$) are of order $O(\ve)$. The
next lemma gives the asymptotic ratio of the widths:

\

\textbf{Lemma 3.1.} \textit{Let $x_0(z)$ and the points
$a(z),b(z),c(z)$ be defined as above. Then}

$$\lim\li_{\ve \da 0}\dfrac{L(a(z),b(z))}{L(b(z),c(z))}=\dfrac{\play{\iint_{D_1(0,z)}\grad\times(\grad F\times b)\cdot
\textbf{n}dm}}{\play{\iint_{D_3(0,z)}\grad\times(\grad F\times
b)\cdot \textbf{n}dm}} \ . \eqno(3.4)$$

\textit{Here the domains $D_1(0,z)$ and $D_3(0,z)$ are the regions
bounded by $C_1(0,z)$ and $C_3(0,z)$.}

\

The \textit{proof} of the lemma is similar to the proof of Lemma 3.4
in [4] but based on (3.3), rather than on the divergency theorem in
$\R^2$, as in [4]. We provide the details in the Appendix.1.
$\square$

\

In the following we will fix an initial point $x$ (not necessarily
$x_0(z)$) on $S(z)$. We put $\widehat{x}=f(x)\in
f(S(z))=R=\{(a,b)\in \R^2: a^2+b^2\leq 1\}$. Let us consider the
trajectory $X_t^\ve(x)$ of (2.7) starting from point $x$. Let
$\widehat{X}_t^\ve(\widehat{x})=f(X_t^\ve(x))$. Our goal now is to
estimate the time of "one rotation" of $X_t^\ve(x)$ around either
$O_1'(z)$ or $O_3'(z)$ or around both of them.

Note that (in two dimensional case), a neighborhood $U$ of a saddle
point of $G$ on $S(z)$ exists such that the system can be reduced to
a linear one in $\widehat{U}\subset \R^2$ by a non-singular
diffeomorphism of the class $C^{1,\al}$, $\al>0$. This comes from
the corresponding result in $\R^2$ ([17, Theorem 7.1]) and the fact
that our surface $S(z)$ is $C^\infty$-diffeomorphic to $R=\{(a,b)\in
\R^2: a^2+b^2\leq 1\}$.

In our case, the system depends on a parameter $\ve$, but one can
check that neighborhood $U$ and $\al>0$ can be chosen the same for
all small enough $\ve$, and the $C^{1,\al}$-norm of the functions
defining the diffeomorphism are bounded uniformly in $\ve$.

For the reason above, it is sufficient to consider the corresponding
flow $\widehat{X}_t^\ve(\widehat{x})$ on $R$. Such a flow has the
same structure consisting of grey and white ribbons on $R$. For
notational convenience we will use the same symbols for objects
related to such a flow, corresponding to our original $X_t^\ve(x)$.
For example, we will write $\widehat{X}_t^\ve(\widehat{x})$ simply
as $X_t^\ve(x)$, and the set $f(U_\dt(x))$ as $U_\dt(x)$, etc. . The
reader could easily understand which specific flow we are referring
to from the context.

The system on $R$ can be linearized in a neighborhood of $O_2'(z)$,
as described above.

First, note that if $x$ is situated outside a fixed (independent of
$\ve$) neighborhood of the $\infty$-shaped curve $C(0,z)$, the
trajectory $X_t^\ve(x)$ comes back to corresponding curve $\gm\ni
x$, orthogonal to the perturbed trajectory , at least, if $\ve>0$ is
small enough. The time of such a rotation $t_{\ve}(x)<\ve A(x)$
(recall that we made time change $t\ra \dfrac{t}{\ve}$); $A(x)$ here
is independent of $\ve$ and bounded uniformly in each compact set
disjoint with $C(0,z)$.

If $x$ is close to $C(0,z)$, then $X_t^\ve(x)$ comes to a
$\dt$-neighborhood $U_\dt(O_2'(z))$ of $O_2'(z)$ in a time less than
$\ve A_\dt$, $A_\dt<\infty$. But the time spent by the trajectory
inside the neighborhood $U_\dt(O_2'(z))$ of $O_2'(z)$ can be large
for small $\ve$; in particular, the separatrices entering $O_2'(z)$
never leave $U_\dt(O_2'(z))$. So we should consider trajectories
started at distance $\dt$ from $O_2'(z)$ in more detail.

Let $\dt>0$ be so small that $U_{2\dt}(O_2'(z))$, for $\ve$ small
enough, belongs to the neighborhood $U$ of $O_2'(z)$ where our
perturbed system can be linearized. The saddle point $O_2'(z)$ under
this transformation goes to the origin $O$, the separatrices of
$O_2'(z)$ go to the axis $\hat{x}$ and $\hat{y}$, the trajectories
$X_t^\ve$ go to the trajectories of the linear system (Fig.4).

One can explicitly calculate the time $\tht(\hat{h},\hat{\dt})$
which the linear system trajectory needs to go from a point
$(\hat{h},\hat{y}_0)$ to $(\hat{\dt}, \hat{y}_1)$ (Fig.4):

$$\tht(\hat{h},\hat{\dt})=\text{const}\cdot
\left|\ln\dfrac{\hat{h}}{\hat{\dt}}\right| \, . \eqno(3.5)$$

\begin{figure}
\centering
\includegraphics[height=6cm, width=6cm , bb=127 32 435 304]{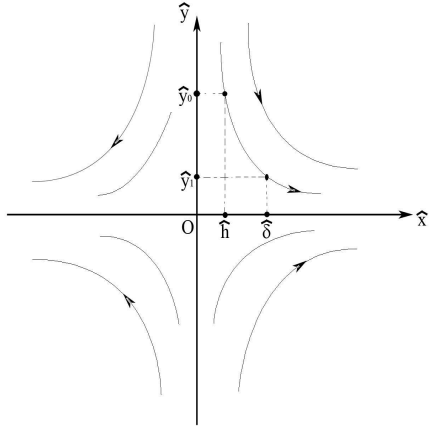}
\caption{Linearized system}
\end{figure}

\

Let a perturbed trajectory enters $U_{\dt}(O_2'(z))$ at a point $x
\in \pt U_{\dt}(O_2'(z))$, $G(x)>0$, and exits $U_\dt(O_2'(z))$ at a
point $y\in \pt U_\dt(O_2'(z))$. We can assume that $x$ and $y$ are
close enough to the pieces of the separatrices which go to the
axises $\hat{x}$, $\hat{y}$ after the linearization so that the
curves $\gm$ and $\gm'$ orthogonal to perturbed trajectories and
containing $x$ and $y$ respectively cross these pieces of
separatrices (these pieces are shown in Fig.5 as bold lines and
denoted by numbers 1,2,3,4) at points $a$ and $a'$ (Fig.5). Let the
distance between $x$ and the closest last piece of the separatrix
entering $O_2'(z)$ be equal to $h$ (here and below we are using the
distance defined by minimal geodesics since we are working in a
sufficiently small neighborhood). Consider the closest to $x$
separatrix crossing $\gm$ at a point $b$ such that $G(b)>G(x)$. Let
$l$ be the distance between $y$ and this separatrix.

\begin{figure}
\centering
\includegraphics[height=8cm, width=8cm , bb=51 3 379 293]{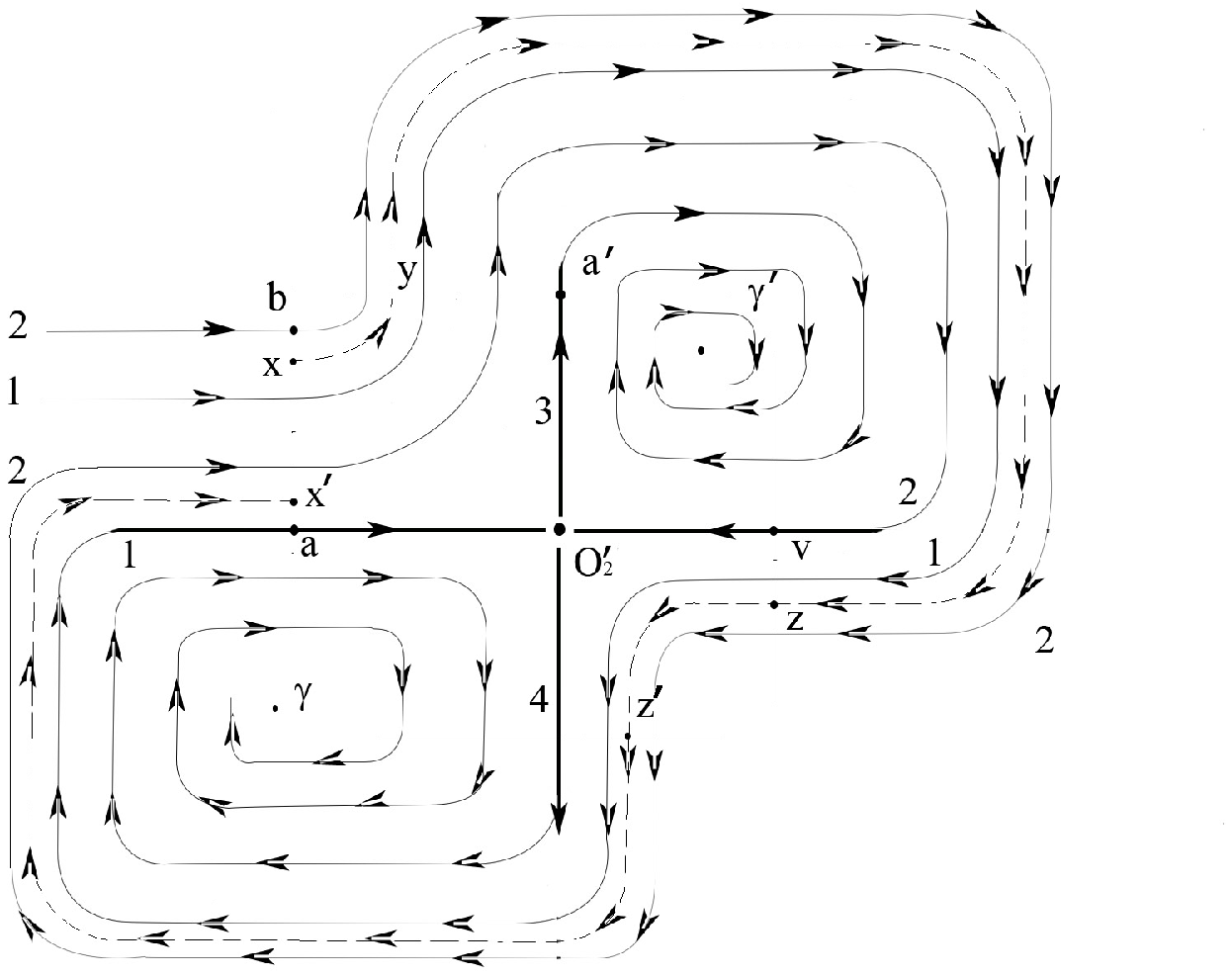}
\caption{Case 1}
\end{figure}

\

If at least one whole ribbon intersects the curve $\gm$ between $x$
and the piece of the separatrix entering $O_2'(z)$ (and containing
point $a$), the trajectory $X_t^\ve(x)$ makes a complete rotation
around both $O_1'(z)$ and $O_3'(z)$ and crosses $\gm$ at a point
$x'\in \gm$ (case 1). The time spent by this trajectory outside
$U_\dt(O_2'(z))$ is bounded from above by $A_1\ve$. Since the
perturbed system can be linearized in $U_{2\dt}(O_2'(z))$ by a
$C^{1,\al}$-diffeomorphism, equality (3.5) implies that the
transition from $x$ to $y$ takes time less than $A_2\ve|\ln h|$;
$A_1$ and $A_2$, in particular, depend on $\dt$, but are independent
of $\ve$.

The trajectory $X_t^\ve(x)$ comes to $\pt U_\dt(O_2'(z))$ again at
the point $z$ (Fig.5). It follows from the divergence theorem that
the distance from $z$ to the last piece of the separatrix entering
$O_2'(z)$ (and containing the point $v$ in Fig.5), in the case when
$X_t^\ve(x)$ comes back to $x'\in \gm$, is bounded from below and
from above by $A_3 h $ and $A_4 h$ respectively. Therefore the
transition from $z$ to $z'$ also takes time less than $A_5 \ve |\ln
h|$.

\begin{figure}
\centering
\includegraphics[height=8cm, width=8cm , bb=47 10 370 288]{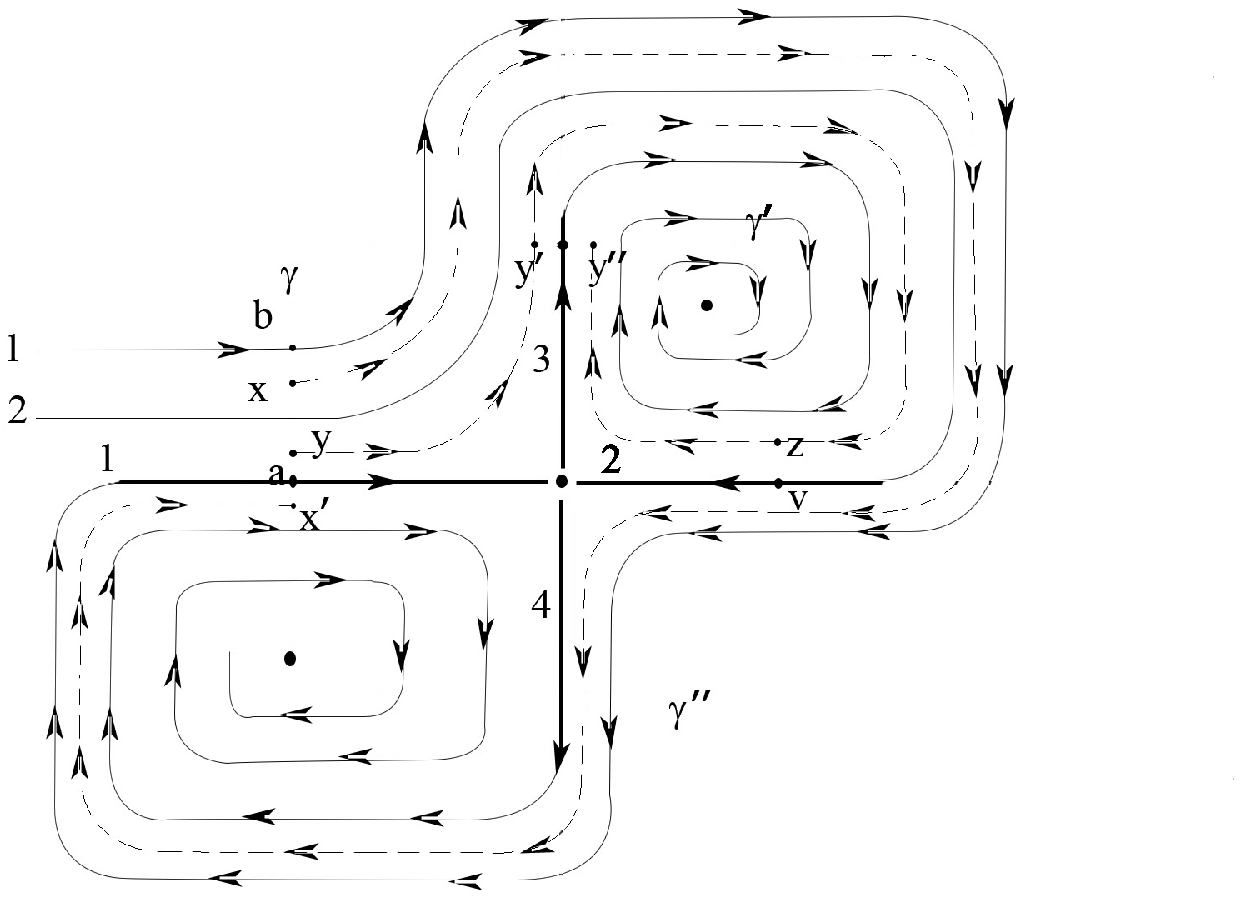}
\caption{Case 2}
\end{figure}

Consider now the case when between the initial point $y\in \pt U_\dt
(O_2'(z))$ and the last piece of the separatrix entering $O_2'(z)$
there is no whole ribbon (Fig.6). Transition between $y$ and $y'$,
because of the same reasons as above, takes time less than
$A_6\ve|\ln h|$, where $h$ is distance between $y$ and the last
piece of separatrix entering $O_2'(z)$. But complete rotation of the
trajectory $X_t^\ve(y)$ includes also the transition from $z$ to
$y''$. It is easy to check using divergence theorem that, the
distance from $z$ to the separatrix entering $O_2'(z)$ is bounded
from below and from above by $A_7 l$ and $A_8 l$ respectively, where
$l$ is the distance between $y$ and the separatrix crossing $\gm$ at
a point $b$ such that $H(b)>H(y)$ (Fig.6). Therefore, the transition
time between $z$ and $y''$ is less than $A_9\ve |\ln l|$, and the
whole rotation time for $X_t^\ve(y)$ is less than $A_{10}\ve(|\ln
h|+|\ln l|)$ for $\ve>0$ small enough.

Denote by $t_{\ve}(x)$ the time of complete rotation for the
trajectory $X_t^\ve(x)$. Suppose $x$ is not a critical point of $G$.
We have

$$t_{\ve}(x)=\min\{t>0: X_t^\ve(x) \text{ crosses twice one of the curves } \gm \text{ or } \gm'\} \, .$$

Summarizing the above bounds and taking into account that outside
$U_\dt(O_1'(z))\cup U_\dt(O_2'(z))\cup U_\dt(O_3'(z))$ the
trajectory $X_t^\ve(x)$ moves with the speed of order $\ve^{-1}$, we
get,

\

\textbf{Lemma 3.2.} \textit{Let $X_t^\ve(x)$ enters $U_\dt(O_2'(z))$
at a point $y=y(x)\in \pt U_\dt(O_2'(z))$, and let $h=h(x)$ be the
distance between $y(x)$ and the last piece of a separatrix entering
$O_2'(z)$. Let $\gm$ be the curve orthogonal to perturbed
trajectories and containing $y(x)$.}

\textit{If in one complete rotation,  $X_t^\ve(y(x))$ come back to
$\gm$, then }
$$t_\ve(x)\leq A_{11}\ve|\ln h(x)| \, . \eqno(3.6)$$

\textit{If $X_t^\ve(y(x))$ does not come back to $\gm$, and $l(x)$
is the distance from $y(x)$ to the closest separatrix, which crosses
$\gm$ at a point $b$, such that $G(b)>G(y(x))$, then for $\ve>0$
small enough, }

$$t_\ve(x)<A_{12}\ve(|\ln h(x)|+|\ln l(x)|) \, . \eqno(3.7)$$

\

Now we come back to our original system (2.7) on $S(z)$. Let $\al$
be a small positive number. Denote by $\cE_{\al}=\cE_\al(\ve)$ the
set of points $x\in S(z)$ such that the distance between $x$ and the
closest separatrix is greater than $\ve\al$ (since $\ve$ is small we
can work with minimal geodesics). Let $\cE_{\al}^{\text{g}}$ be the
intersection of $\cE_\al$ with the gray ribbon; $\cE_\al^{\text{w}}$
be the intersection with the white ribbon.

Denote by $\Lb_\ve(x,\bt)$ the time when $X_t^\ve(x)$ reaches
$C(\bt,z)$:

$$\Lb_\ve(x,\bt)=\inf\{t>0: G(X_t^\ve(x))=\bt\} \, ;$$
if $G(x)>0$ and $|\bt|$ is small, $\Lb_\ve(x,\bt)<\infty$ for all
small $\ve>0$.

\

\textbf{Lemma 3.3.} \textit{Let $G(x_0(z))>0$ and let $\mu>0$ be so
small that $G(x)>0$ for $x\in U_{2\mu}(x_0(z))$. There exist
$\al_0$, $\bt_0>0$ and $A_{13}$ such that for each $x\in
U_{\mu}(x_0(z))\cap\cE_\al$, $\al\in (0,\al_0)$, $\bt\in (0,\bt_0)$,
}
$$\Lb_\ve(x,-\bt)-\Lb_\ve(x,\bt)<A_{13}\bt|\ln \bt|  \eqno(3.8)$$
\textit{for $\ve<\ve_0$. Here $A_{13}$, in particular, depends on
$\al$ and $\bt$ but is independent of $\ve$; $\ve_0>0$ depends on
$\al$ and $\bt$.}

\

The \textit{proof} of this lemma is based on Lemma 3.2 and the fact
that each rotation decreases the value of $G$ on an amount of order
$O(\ve)$. Therefore the total time is less than
$A_{14}\sum\li_{k=1}^{[\dfrac{\bt}{\ve}]}\ve|\ln(k\ve)|\sim
\play{\int_0^\bt |\ln z| dz\leq A_{13}\bt |\ln \bt|}$ for $\ve>0$
small enough. $\square$

\

\textit{Proof of Theorem 3.1.} Equation (3.2) can be considered for
each of three edges of the graph $\Gm$ corresponding to $G(x)$ on
$S(z)$: for $i=1,2,3$, we have

$$\begin{array}{l}
\dot{g}_t^{(i)}=\dfrac{1}{T_i(g_t^{(i)})}B^{(i)}(g_t^{(i)},z) \, ,
\\  \\ \play{T_i(g)=\oint_{C_i(g,z)}\dfrac{dl}{|\grad F\times \grad G|}} \,
, \\
\play{B^{(i)}(g,z)=-\iint_{D_i(g,z)}\grad \times (\grad F \times
\mathbf{b})\cdot \mathbf{n} dm \ .}
\end{array} \eqno(3.10)$$

Equation (3.10) for $i=2$ can be solved for each initial condition
$g_0^{(2)}=g>0$, $g<\max\{G(w): w\in \pt S(z)\}$. Such a solution is
unique, and $g_t^{(2)}$ reaches $0$ in a finite time $\tau_0(g,z)$.
If $i=1,3$, equation (3.10) with initial condition $g_0^{(i)}=g<0$
has a unique solution; if $g_0^{(i)}=0$, equation (3.10) has a
unique solution $\widetilde{g}_t^{(i)}$ if we additionally assume
that $\widetilde{g}_t^{(i)}<0$ for $t>0$.

Define two continuous functions $\widehat{g}_t^1(g)$ and
$\widehat{g}_t^3(g)$, $t\geq 0$, as follows:
$\widehat{g}_0^1=\widehat{g}_0^3=g>0$,

$$\widehat{g}_t^1(g)=\left\{
\begin{array}{l}
g_t^{(2)} \, , \, g_0^{(2)}=g \, , \, 0\leq t \leq \tau_0(g,z) \, ,
\\ \widetilde{g}_{t-\tau_0(y,z)}^{(1)} \, , \, \tau_0(g,z)\leq t < \infty  \, ;
\end{array}
\right.
$$

$$\widehat{g}_t^3(g)=\left\{
\begin{array}{l}
g_t^{(2)} \, , \, g_0^{(2)}=g \, , \, 0\leq t \leq \tau_0(g,z) \, ,
\\ \widetilde{g}_{t-\tau_0(g,z)}^{(3)} \, ,  \, \tau_0(g,z)\leq t < \infty  \, ;
\end{array}
\right.
$$

Let us cut out $\al\ve$-neighborhoods of the separatrices
($\mu$-neighborhood of a point $x_0$, $G(x_0)>0$, is shown in
Fig.7); recall that $\cE_\al$ is the exterior of the
$\ve\al$-neighborhood of the separatrices, $\cE_\al^\text{g}$ is the
intersection of $\cE_\al$ with the gray ribbon, $\cE_\al^\text{w}$
is the intersection of $\cE_\al$ with the white ribbon. In
particular, $\cE_0^\text{g}$ ($\cE_0^\text{w}$) is whole gray
(white).

\begin{figure}
\renewcommand{\captionlabeldelim}{.}
\centering
\includegraphics[height=8cm, width=9cm , bb=80 23 347 254]{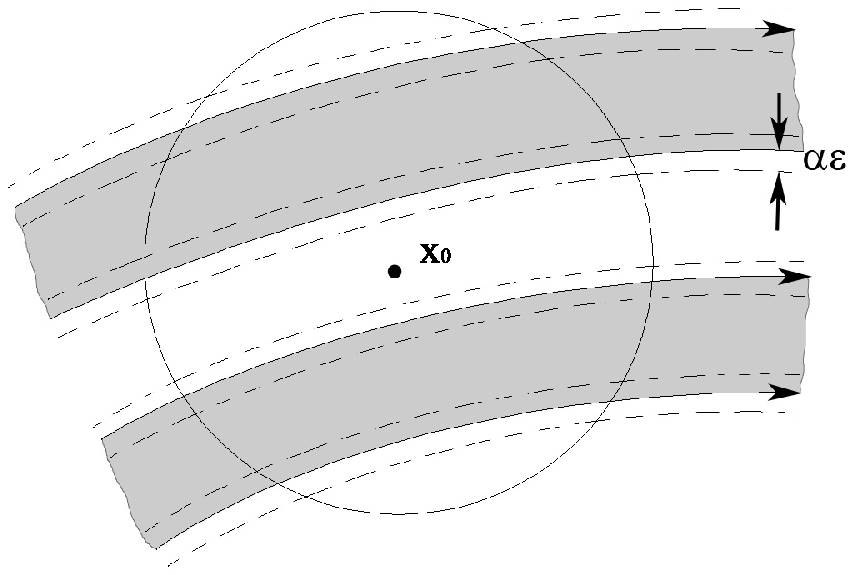}
\caption{}
\end{figure}

The classical averaging principle together with Lemma 3.3 imply that
for each $x\in U_\dt(x_0(z))\cap \cE_\al^\text{g}$, $G(x)=g>0$, for
any $\lb, T>0$, and any small enough $\al, \dt>0$, there exists
$\ve_0>0$ such that

$$\max\li_{0\leq t \leq T}|H(X_t^\ve(x))-\widehat{g}_t^1(g)|<\lb \eqno(3.11)$$
for $0<\ve<\ve_0$.

Similarly, for each $x\in U_\dt(x_0(z))\cap \cE_\al^\text{w}$,
$G(x)=g>0$,

$$\max\li_{0\leq t \leq T}|H(X_t^\ve(x))-\widehat{g}_t^3(g)|<\lb \eqno(3.12)$$
for $0<\ve<\ve_0$.

Let $G(x)>0$ for $x\in U_\dt(x_0(z))$ so that
$\iY(U_\dt(x_0(z)))\subset I_2 \subset \Gm$. Define a stochastic
process $Y_t^\dt(x_0(z))$, $t\geq 0$, on $\Gm$ as follows: (recall
that the pair $(k,g)$, where $k$ is the number of an edge, $k\in
\{1,2,3\}$, and $g$ is the value of $G(x)$ on $\iY^{-1}(y)$, $y\in
\Gm$, form a global coordinate system on $\Gm$)

$$Y_t^\dt(x_0(z))=(2, \widehat{g}_t^2(G(x_0(z,\dt)))) \text{ for } 0\leq t \leq \tau_0(x_0(z,\dt)) \, .$$
(Recall that $\tau_0(x_0(z,\dt))$ is the first time when the process
$g_t^{(2)}$, $g_0^{(2)}=G(x_0(z,\dt))>0$ in (3.10) reaches $0$.)

At the time $\tau_0(x_0(z,\dt))$ the process $Y_t^\dt(x_0(z,\dt))$
reaches $O_2(z)$ and without any delay goes to $I_1$ or $I_3$ with
probabilities
$$p_1=\dfrac{\play{\iint_{D_1(0,z)}\grad\times(\grad F\times b)\cdot
\textbf{n}dm}}{\play{\iint_{D_1(0,z)}\grad\times(\grad F\times
b)\cdot \textbf{n}dm}+\play{\iint_{D_3(0,z)}\grad\times(\grad
F\times b)\cdot \textbf{n}dm}} \ , \eqno(3.13)
$$
$$p_3=\dfrac{\play{\iint_{D_3(0,z)}\grad\times(\grad F\times
b)\cdot \textbf{n}dm}}{\play{\iint_{D_1(0,z)}\grad\times(\grad
F\times b)\cdot
\textbf{n}dm}+\play{\iint_{D_3(0,z)}\grad\times(\grad F\times
b)\cdot \textbf{n}dm}} \ , \eqno(3.14)$$ respectively;
$Y_t^\dt(x_0(z))=(1,
\widehat{g}^1_{t-\tau_0(x_0(z,\dt))}(x_0(z,\dt))$ for
$\tau_0(x_0(z,\dt))\leq t < \infty$ if $Y_t^\dt(x_0(z))$ enters
$I_1$ at time $\tau_0(x_0(z,\dt))$, and $Y_t^\dt(x_0(z))=(3,
\widehat{g}^3_{t-\tau_0(x_0(z,\dt))}(x_0(z,\dt))$ for
$\tau_0(x_0(z,\dt))\leq t < \infty$ if $Y_t^\dt(x_0(z,\dt))$ enters
$I_3$ at time $\tau_0(x_0(z,\mu))$.

One can consider a process $\overline{Y}_t(x_0(z))=Y_t^0(x_0(z))$ on
$\Gm$: $\overline{Y}_t(x_0(z))$ is deterministic inside the edges
and governed by equations (3.10); its stochasticity concentrated at
the vertex $O_2(z)$: after reaching $O_2(z)$,
$\overline{Y}_t(x_0(z))$ immediately goes to $I_1$ or to $I_3$ with
probabilities $p_1$ or $p_3$ defined by equalities (3.13) and
(3.14).

Denote by $\text{Area}(D)$, $D\subset S(z)$, the area of a domain
$D$. Since the point $x_0(z,\dt)$ is distributed uniformly in
$U_\dt(x_0(z))$,

$$\begin{array}{l}
\play{\left|\Prob\{X_t^{\ve,\dt} \text{ enters }
D_1(0,z)\}-\dfrac{\text{Area}(\cE_0^\text{g}\cap
U_\dt(x_0(z)))}{\text{Area}(U_\dt(x_0(z)))}\right|\ra 0 \, ,}
\\
\\
\play{\left|\Prob\{X_t^{\ve,\dt} \text{ enters }
D_3(0,z)\}-\dfrac{\text{Area}(\cE_0^\text{w}\cap
U_\dt(x_0(z)))}{\text{Area}(U_\dt(x_0(z)))}\right|\ra 0 \, ,}
\end{array} \eqno(3.15)$$
as $\ve \da 0$. According to Lemma 3.2,

$$\lim\li_{\ve \da 0}\dfrac{\text{Area}(\cE_0^\text{g}\cap U_\mu(x_0))}{\text{Area}(U_\mu(x_0))}=p_1
\, , \, \lim\li_{\ve \da 0}\dfrac{\text{Area}(\cE_0^\text{w}\cap
U_\mu(x_0))}{\text{Area}(U_\mu(x_0))}=p_3 \, , \eqno(3.16)$$ where
$p_1$ and $p_3$ are defined in (3.13) and (3.14).

Taking into account that $\text{Area}(\cE_\al^\text{g}\cap
U_\dt(x_0(z)))\ra \text{Area}(\cE_0^\text{g}\cap U_\dt(x_0(z)))$ and
$\text{Area}(\cE_\al^\text{w}\cap U_\dt(x_0(z)))\ra
\text{Area}(\cE_0^\text{w}\cap U_\dt(x_0(z)))$ as $\al \da 0$, we
derive from (3.10)-(3.16) that, for each $T>0$, the slow component
$\iY(X_t^{\ve,\dt})$ of $X_t^{\ve,\dt}$ converges weakly in the
space of continuous functions on $[0,T]$ with values in $\Gm$ to the
process $Y_t^\dt(x_0(z))$.

It is easy to see that $Y_t^\dt(x_0(z))$ converges weakly to
$\overline{Y}_t(x_0(z))$ as $\dt \da 0$.

This gives the proof of Theorem 3.1. $\square$

\

\section{Regularization by stochastic perturbation of the dynamics}

Let now the perturbation have deterministic and stochastic parts:

$$\dot{X}_t^{\ve,\dt}=\dfrac{1}{\ve}(\grad F \times \grad G)(X_t^{\ve,\dt})+
\grad F \times (\grad F \times
\mathbf{b})(X_t^{\ve,\dt})+\dt\sm(X_t^{\ve,\dt})\circ \dot{W}_t \ ,
\eqno(4.1)$$ $X_0^{\ve,\dt}=x_0(z)\in \R^3$, $F(x_0(z))=z$, and
$\dt>0$, $0<\ve<<1$. The stochastic term $\sm(X_t^{\ve,\dt})\circ
\dot{W}_t$ is understood in the Stratonovich sense. The $3\times 3$
matrix $\sm(x)=(\sm_{ij}(x))$ is assumed to be smooth and satisfy
the relation $\sm^T\grad F\equiv0$. If we denote by
$a(x)=(a_{ij}(x))=\sm(x)\sm^T(x)$ the diffusion matrix, the
condition $\sm^T(x)\grad F (x)=0$ is equivalent to the assumption
that $a(x)\grad F (x)\equiv 0$. By using the It\^{o} formula for
Stratonovich integrals, we have, that
$$\dfrac{dF(X_t^{\ve,\dt})}{dt}=\grad F \cdot
\left[\dfrac{1}{\ve}(\grad F \times \grad G)(X_t^{\ve,\dt})+ \grad F
\times (\grad F \times
\mathbf{b})(X_t^{\ve,\dt})+\dt\sm(X_t^{\ve,\dt})\circ
\dot{W}_t\right]=0 \ .$$ (One can directly check that equality
$\sm^T\grad F=0$ implies $\grad F \cdot \sm\circ \dot{W}_t=0 \ .$)
In particular, if $\mathbf{b}(x)\equiv 0$, we have pure stochastic
perturbations. Therefore $F$ is a first integral for system (4.1),
i.e., the process $X_t^{\ve,\dt}$ never leaves the surface
$\{F=z\}$. We also assume that $\mathbf{e}\cdot (a(x)\mathbf{e})\geq
\underline{a} |\mathbf{e}|^2$ for a constant $\underline{a}>0$ and
every $\mathbf{e}\in \R^3$ such that $\mathbf{e}\cdot\grad F (x)=0$.
This means that the process $X_t^{\ve,\dt}$ is non-degenerate if
considered on the manifold $S(z)\subset \{x\in \R^3: F(x)=z\}$.
Recall that we stop our process $X_t^{\ve,\dt}$ once it hits $\pt
S(z)$. The resulting process is still called $X_t^{\ve,\dt}$.

We will make use of the following simple Lemma (see, for instance,
[20, page 36, formula (3.3.6)]):

\

\textbf{Lemma 4.1.} \textit{Let $\mathbf{b}(s,x): \R_+\times \R^d\ra
\R^d$ be Lipschitz and bounded in $s$ and $x$. Let $\sm(s,x) :
\R_+\times \R^d \ra \R^d\times \R^d$ be bounded, Lipschitz in s, and
differentiable in $x$. Let $\sm_{ij}$ is the $(i,j)$-th element of
matrix $\sm$. Consider the diffusion process
$$X_t=x+\int_0^t\mathbf{b}(s,X_s)ds+\int_0^t\sm(s,X_s)\circ dW_s $$ in $\R^d$, where the
stochastic term is understood in Stratonovich sense. Then we have
$$X_t=x+\int_0^t \mathbf{b}(s,X_s)ds+\dfrac{1}{2}\int_0^t \mathbf{c}(s,X_s)ds+\int_0^t
\sm(s,X_s)dW_s$$ where the stochastic term is understood in the
It\^{o} sense. Here vector $\mathbf{c}(s,x)\in \R^d$ has $i$-th
component $\sum\li_{j,k=1}^d\dfrac{\pt \sm_{ij}}{\pt x_k}\sm_{kj}$,
$1\leq i \leq d$.}

\

Using this Lemma, we easily write equation (4.1) in the It\^{o}
sense:
$$\dot{X}_t^{\ve,\dt}=\dfrac{1}{\ve}(\grad F \times \grad G)(X_t^{\ve,\dt})+
\grad F \times (\grad F \times
\mathbf{b})(X_t^{\ve,\dt})+\dt\sm(X_t^{\ve,\dt})\dot{W}_t +
\dfrac{\dt^2}{2}\mathbf{\Sm}(X_t^{\ve,\dt})\ . \eqno(4.2)$$

Here $\mathbf{\Sm}$ is a vector in $\R^3$ with the $i$-th component
$\Sm_i=\sum\li_{j,k=1}^3\dfrac{\pt \sm_{ij}}{\pt x_k}\sm_{kj}$ for
$i=1,2,3$.

The generator $L$ of the process $X_t^{\ve,\dt}$ is written as

$$Lu(x)=\dfrac{\dt^2}{2}\sum\li_{i,j=1}^{3}a_{ij}(x)\dfrac{\pt^2
u}{\pt x_i\pt x_j}+\left(\dfrac{1}{\ve}\grad F \times \grad G +\grad
F \times (\grad F \times
\mathbf{b})+\dfrac{\dt^2}{2}\mathbf{\Sm}\right)\cdot \grad u(x) \ .
\eqno(4.3)$$

Using It\^{o}'s formula we see that

$$\begin{array}{l}
\play{G(X_t^{\ve,\dt})-G(x_0(z))=\dt\int_0^t(\grad
G)^T(X_s^{\ve,\dt})\sm(X_s^{\ve,\dt}) dW_s+}\\\play{\ \ \ \ \ \ \ \
\ \ \ +\int_0^t\left(\grad G \cdot (\grad F \times (\grad F \times
\mathbf{b}))+ \dfrac{\dt^2}{2}\grad G \cdot \mathbf{\Sm }+
\dfrac{\dt^2}{2}\sum\li_{i,j=1}^3a_{ij}\dfrac{\pt^2 G}{\pt x_i \pt
x_j}\right)(X_s^{\ve,\dt})ds} .\end{array} \eqno(4.4)$$

Now we are in a position to use the standard averaging principle
(see, for example, [14, Chapter 8]), to check that within edge $I_i$
($i=1,2,3$) of the graph $\Gm$, as $\ve \da 0$ and $\dt$ is fixed,
the process $G(X_t^{\ve,\dt})$ converges weakly to the process
$G_t^\dt$ governed by the operator

$$\overline{L}_i=\dfrac{1}{T_i(g)}\left(A_i(g,z)+\dfrac{\dt^2}{2}A_{1,i}(g,z)+
\dfrac{\dt^2}{2}A_{2,i}(g,z)\right)\dfrac{d}{dg}+ \dfrac{\dt^2}{2}
\dfrac{1}{T_i(g)}B_i(g,z)\dfrac{d^2}{dg^2} \ . \eqno(4.5)$$

The coefficients are

$$\begin{array}{l}
\play{A_i(g,z)=\oint_{C_i(g,z)}\grad G \cdot (\grad F \times (\grad
F\times \mathbf{b}))\dfrac{dl}{|\grad F \times \grad G|} \ ,}\\
\play{A_{1,i}(g,z)=\oint_{C_i(g,z)}\grad G \cdot \mathbf{\Sm}
\dfrac{dl}{|\grad F \times \grad G|} \ ,}\\
\play{A_{2,i}(g,z)=\oint_{C_i(g,z)}\sum\li_{i,j=1}^3a_{ij}\dfrac{\pt^2G}{\pt
x_i\pt x_j} \dfrac{dl}{|\grad F \times \grad G|} \ ,}\\
\play{B_i(g,z)=\oint_{C_i(g,z)}|(\grad G)^T\sm|^2\dfrac{dl}{|\grad F
\times \grad G|} \ .}\end{array} \eqno(4.6) $$

Here $T_i(g)$ is the period of rotation of the unperturbed system
(2.1) on $C_i(g,z)$: $\play{T_i(g)=\oint_{C_i(g,z)}\dfrac{dl}{|\grad
F \times \grad G|}}$, where $dl$ is the length element on
$C_i(g,z)$.

We define a process $Y_t^\dt$ on $\Gm$ as follows: $Y_t^\dt$ is a
Markov process on $\Gm$, stopped once it hits exterior vertex $P$
(recall that we stop our process $X_t^{\ve,\dt}$ once it hits $\pt
S(z)$; also recall that by our definition $S(z)=\{x\in \R^3:
G(x)\leq G(x_0(z))+1\}\cap \{x\in \R^3: F(x)=z\}$ so that $\pt S(z)$
is a trajectory, corresponding to vertex $P$ on $\Gm$) and governed
by a generator $A$. The operator $A$ is defined as follows. The
domain of definition for the generator $A$ consists of functions
$f(g,k)$ on $\Gm$ which are twice continuously differentiable in the
variable $g$ within the interior part of each edge $I_i$; inside
$I_i$, $Af(g,i)=\overline{L_i} f(g,i)$, and finite limits
$\lim\li_{y\ra O_i(z)}Af(y)$ (which are taken as the value of $Af$
at vertex $O_i(z)$) and finite one sided limits $\lim\li_{g\ra
G(O_i(z))}\dfrac{\pt f}{\pt g}(g,i)$, $\lim\li_{g\ra G(P)}\dfrac{\pt
f}{\pt g}(g,i)$ exist. We set $\lim\li_{y\ra P}Af(y)=0$ (taken as
the value of $Af$ at point $P$, this means that the process
$Y_t^\dt$ is stopped at the point $P$). For the interior vertex
$O_2(z)$, $f$ satisfies the gluing condition:
$$\sum\li_{i=1}^3(\pm)\bt_{2,i}\lim\li_{g\ra G(O_2(z))}\dfrac{\pt f}{\pt
g}(g,i)=0 \ , \eqno(4.7)$$ where $+$ sign is for the limit taking
within edge $I_2$ and $-$ sign is for the limit taking within edge
$I_1$ and $I_3$. The coefficients $\bt_{2,i}$ are defined by
$$\bt_{2,i}=\oint_{C_i(0,z)}|(\grad G)^T\sm|^2\dfrac{dl}{|\grad F
\times \grad G|} \ . \eqno(4.8)$$ Exterior vertex $O_1(z)$ and
$O_3(z)$ are inaccessible. Such a process $Y_t^\dt$ on $\Gm$ exists
and is unique ([14, Chapter 8]).

\

\textbf{Theorem 4.1.} \textit{As $\ve \da 0$ and $\dt$ is fixed, the
process $\iY(X_t^{\ve,\dt})$ converges weakly in the space of
continuous functions $f: [0,T]\ra \Gm$, $0<T<\infty$, to the process
$Y_t^\dt$.}

\

The \textit{proof} of this Theorem is based on the fact that we can
carry the dynamics of (3.1) on $S(z)$ to a corresponding one on
$R\subset\R^2$ by the $C^\infty$-diffeomorphism $f:S(z)\ra \R^2$. We
denote $f(x_1,x_2,x_3)=(f_1(x_1,x_2,x_3),f_2(x_1,x_2,x_3))$ ,
$(x_1,x_2,x_3)\in S(z)$. Let $Z_t^{\ve,\dt}=f(X_t^{\ve,\dt})$ be the
image of the diffusion process on $\R^2$. Using the It\^{o} formula
for Stratonovich integrals, we have
$$\begin{array}{l}
dZ_t^{\ve,\dt} \\
=df(X_t^{\ve,\dt})
\\=(Df)(f^{-1}(Z_t^{\ve,\dt}))dX_t^{\ve,\dt}
\\=\left(\dfrac{1}{\ve}\vec{\bt}(Z_t^{\ve,\dt})+\vec{\bt_1}(Z_t^{\ve,\dt})\right)dt+\dt (Df)(f^{-1}(Z_t^{\ve,\dt}))
\sm(f^{-1}(Z_t^{\ve,\dt}))\circ dW_t
\\=\left(\dfrac{1}{\ve}\vec{\bt}(Z_t^{\ve,\dt})+\vec{\bt_1}(Z_t^{\ve,\dt})\right)dt+
\dt \widetilde{\sm}(Z_t^{\ve,\dt})\circ d\widetilde{W}_t \,
\end{array} \eqno(4.9)$$ so that $Z_t^{\ve,\dt}=(f_1(X_t^{\ve,\dt}),f_2(X_t^{\ve,\dt}))$ is a
diffusion process on $\R^2$, stopped once it hits $\pt R$.

Here the matrix $Df$ is the differential of $f$:
$Df=\left(\dfrac{\pt f_i}{\pt x_j}\right)_{1\leq i\leq 2, 1\leq
j\leq 3}$. The vector fields $\vec{\bt}(Z)=(Df)(\grad F \times \grad
G)(f^{-1}(Z))$ and $\vec{\bt}_1(Z)=(Df)(\grad F \times (\grad F
\times \mathbf{b}))(f^{-1}(Z))$. The $2\times 2$ matrix
$\widetilde{\sm}$ is defined in the following way:
$\widetilde{\sm}(Z)\circ
d\widetilde{W_t}=(Df)(f^{-1}(Z))\sm(f^{-1}(Z))\circ dW_t$, where
$W_t$ is the standard $3$-dimensional Wiener process and
$\widetilde{W}_t$ is the standard $2$-dimensional Wiener process.
The integral curves of the vector field $\vec{\bt}$ has one saddle
point $f(O_2(z))$ and two stable equilibriums $f(O_1(z))$ and
$f(O_3(z))$.

We define $\cG(Z)=G(f^{-1}(Z))$ for $Z\in \R^2$. The function $\cG$
serves as the first integral for the vector field $\vec{\bt}$:
$\grad \cG \cdot \vec{\bt}=0$. Furthermore, it is easy to check that
$\vec{\bt}(Z)=\kp(Z)\bar{\grad} \cG(Z)$ with $\kp \neq 0$, so that
our system just by a non-singular time change differs from a
Hamiltonian system with one degree of freedom. Therefore one can use
the same arguments as in the case of $2$-dimensional Hamiltonian
systems (see, [14, Chapter 8], [13], [11]) to calculate the limiting
behavior the process $Z_t^{\ve,\dt}$ as $\ve\da 0$. (In the
calculation of the gluing conditions, the problem caused by
additional drift term $\vec{\bt_1}$ and another drift term related
to the Stratonovich integral can be resolved using the absolute
continuous transformation; detailed estimates see [11] and
Appendix.2.) The coefficients of the gluing condition at the
interior vertex are given as follows:

$$\bt_{2,i}=\oint_{f(C_i(0,z))}|(\grad \cG)^T\widetilde{\sm}|^2\dfrac{dl_z}{|\vec{\bt}|} \ ,$$
where $dl_z$ is the length element on $f(C_i(0,z))$. Note that they
coincide with (4.8), since equality

$$(\grad \cG)^T\widetilde{\sm}\circ d\widetilde{W_t}=(\grad G)^T (Df^{-1}) (Df) \sm\circ dW_t=(\grad G)^T \sm\circ dW_t$$
 implies $$|(\grad \cG)^T\widetilde{\sm}|^2=|(\grad G)^T\sm|^2 \ ,$$ and
$\play{\dfrac{dl_z}{|\vec{\bt}|}=\dfrac{dl}{|\grad F \times \grad
G|}} \ .$ $\square$

\

The next step is to consider the limit as $\dt \da 0$ of the process
$Y_t^\dt$. This follows the same line of argument as in [4, Section
2]. In particular, one can do a similar calculation as in Lemma 2.2
of [4]. The additional small drift term depending on $\dt$ (caused
by the Stratonovich integral) in (4.5) will disappear as $\dt \da
0$. (We briefly indicate how to calculate this in Appendix.3.) We
therefore have a limiting process $Y_t$ on $\Gm$ defined as follows:
$Y_t=(g^{(i)}_t,k_t)$ is a deterministic motion inside each edge of
$\Gm$ with $g^{(i)}_t$ satisfying the differential equation (3.10)
and the branching probability for $Y_t$ at vertex $O_2(z)$ is given
by (3.13) and (3.14). The process $Y_t$ spends time zero at the
vertex $O_2(z)$. These arguments imply

\

\textbf{Theorem 4.2.} \textit{As $\dt \da 0$, the process $Y_t^\dt$
converges weakly in the space of continuous functions $f: [0,T] \ra
\Gm$, $0<T<\infty$, to the process $Y_t$.}

\

Theorem 4.1 and 4.2 imply that the slow component
$\iY(X_t^{\ve,\dt})$ of the process $X_t^{\ve,\dt}$ converges weakly
to the process $Y_t$ on the graph $\Gm$. Note that $Y_t$ is
independent of the diffusion matrix $a(x)=\sm(x)\sm^T(x)$ and is the
same process which we had using regularization by stochastic
perturbation of the initial point.

\

Consider process $X_t^{\ve,\dt}$ defined by (4.1). Under the
assumption that the deterministic perturbation in (4.1) is
friction-like, for $\ve>0$ small enough and fixed, the equilibrium
$O_1'(z)$ and $O_3'(z)$ are asymptotically stable for the dynamical
system $X_t^{\ve,0}$ on $\{F(x)=z\}$. The process $X_t^{\ve,\dt}$ is
close to $X_t^{\ve,0}$ on any fixed time interval if $\dt$ is small
enough. But on time intervals of order
$\exp\left\{\dfrac{\lb}{\dt^2}\right\}$ for $\lb>0$, $X_t^{\ve,\dt}$
may perform transitions between the neighborhoods of $O_1'(z)$ and
$O_3'(z)$ due to the large deviations from $X_t^{\ve,0}$. In a
generic case, for $x\in \{F(x)=z\}$ and $\lb>0$, there exists just
one stable equilibrium $M^\ve(x,\lb)$ (in the case of two stable
equilibriums, $M^\ve(x,\lb)=O_1'(z)$ or $M^\ve(x,\lb)=O_3'(z)$) such
that with probability close to 1 as $\dt \da 0$,
$X_{T^{\dt}(\lb)}^{\ve,\dt}$ is situated in a small neighborhood of
$M^\ve(x,\lb)$, if $X_0^{\ve,\dt}=x$, $\lim\li_{\dt \da 0}\dt^2 \ln
T^\dt(\lb)= \lb$. The state $M^\ve(x,\lb)$ is called metastable
state for a given initial point $x$ and time scale $\lb>0$ (see [8],
[10] where the procedure for calculating $M^\ve(x,\lb)$ is
described).

But it turns out that the function $M^\ve(x,\lb)$ is very sensitive
to $\ve$ as $\ve \da 0$: For $\lb$ not very large, $M^\ve(x,\lb)$
alternatively is equal to $O_1'(z)$ or to $O_3'(z)$ as $\ve \da 0$.
Moreover, for small $\ve$, $M^{\ve}(x,\lb)$ is sensitive to changes
of the initial point $x$ as well. Therefore, if $\ve<<1$, the notion
of metastability should be modified (compare with [3], [9]): For
given $x$ and $\lb$, one should consider the set of metastable
distributions between the stable equilibriums. In general, there
exists a finite number of distributions on the set of stable
equilibriums which serve as limiting distributions of
$X_{T^\dt(\lb)}^{\ve,\dt}$ as first $\ve \da 0$ and then $\dt \da
0$. The set of metastable distributions is independent of the
stochastic terms in (4.1) and defined just by the deterministic
system and deterministic perturbations. But which of those
distributions serves as limiting distribution of
$X_{T^\dt(\lb)}^{\ve,\dt}$, $X_0^{\ve,\dt}=x$, is defined by the
stochastic term in (4.1).

In our case, when we have just two stable equilibriums $O_1(z)$ and
$O_3(z)$, three distributions can serve as metastable distribution:
first, the distribution concentrated at $O_1(z)$, second, the
distribution concentrated at $O_3(z)$, and third, the distribution
between $O_1(z)$ and $O_3(z)$ with $\Prob\{O_1(z)\}=p_1$,
$\Prob\{O_3(z)\}=p_3$ where $p_1$ and $p_3$ defined by (3.13),
(3.14).

\

\textbf{Theorem 4.3.} \textit{Let
$\play{\lb_1=-\int_{G(O_1(z))}^0\dfrac{A_1(g,z)dg}{B_1(g,z)}<\int_{G(O_3(z))}^0
\dfrac{A_3(g,z)dg}{B_3(g,z)}}=\lb_3$, where $A_i(g,z)$ and
$B_i(g,z)$ are defined by (4.6). Let $\lim\li_{\dt\da 0}\dt^2\ln
T^{\dt}(\lb)=\lb>0$. Then for each small enough $h>0$, }

$$\lim\li_{\dt \da 0}\lim\li_{\ve \da 0}\Prob_x\{|X_{T^\dt(\lb)}^{\ve,\dt}-O_1(z)|<h\}=1
\text{ \textit{if} } \iY(x)\in I_1 \text{ \textit{and} } \lb<\lb_1
,$$

$$\lim\li_{\dt \da 0}\lim\li_{\ve \da 0}\Prob_x\{|X_{T^\dt(\lb)}^{\ve,\dt}-O_3(z)|<h\}=1
\text{ \textit{if} } \iY(x)\in I_2 \text{ and } \lb>0 \text{
\textit{or if} } \lb>\lb_3 \text{ \textit{for any} } x\in \{F(x)=z\}
, $$

$$\lim\li_{\dt \da 0}\lim\li_{\ve \da 0}\Prob_x\{|X_{T^\dt(\lb)}^{\ve,\dt}-O_i(z)|<h\}=p_i \ , i\in
\{1,3\} \ ,
  \text{ \textit{if} } \iY(x)\in I_3 \text{ \textit{and} } \lb<\lb_1 .$$

\textit{The probabilities $p_1$ and $p_3$ are defined by
(3.13)-(3.14). }

\

The proof follows from Theorem 4.1 and the fact that the transition
time from $O_1(z)$ to $O_3(z)$ (from $O_3(z)$ to $O_1(z)$) for the
process $Y_t^\dt$ on $\Gm$ is logarithmically equivalent as $\dt \da
0$ to $\exp\left\{\dfrac{\lb_1}{\dt^2}\right\}$
($\exp\left\{\dfrac{\lb_3}{\dt^2}\right\}$) (Theorem 4.4.2 in [14]).
$\square$

\

\textbf{Remark:} We assumed in Sections 3 and 4 that the function
$G(x)$ has in $S(z)$ just one saddle point and two minima. We also
assumed that the deterministic perturbations are friction-like. Then
each minimum point become asymptotically stable for the perturbed
system. It is not difficult to check that if $G(x)$ has on the set
$\widetilde{S}(z)=\{F=z\}$ (we assumed it has only one connected
component) more than two minima points and several saddle points but
just one local maximum, and the deterministic perturbations are
friction-like, then the system can be regularized by an addition of
stochastic perturbations of the initial point or of the dynamics.
Corresponding graph in this case has several interior vertices
corresponding to the saddle points of $G(x)$ and exterior vertices
corresponding to the extremums.

Inside each edge, the limiting slow motion is governed by
corresponding equation (3.2). The exterior vertices are inaccessible
in finite time. The limiting slow motion spends time zero at
interior vertices, and the branching at each interior vertex occurs
exactly as in the case of a unique saddle point. The branching at
each interior point is independent of the previous behavior of the
limiting slow motion.

But situation is a bit different if $G(x)$ has on $\widetilde{S}(z)$
more than one maxima or if the perturbations are not friction-like.
In this case, in general, it is impossible to regularize the problem
by a random perturbation of the initial point: the limit of
$\iY(X_t^{\ve,\dt})$ as $\ve \da 0$ may not exist (compare with
[4]). The regularization by stochastic perturbations of the
equation, as we did in Section 4, is possible under mild additional
assumptions. One should keep in mind that, if the deterministic
perturbation is not friction-like, the stochastic branching occurs
just at those interior vertices where there are two "exit" edges and
one "entrance" edge (this means that the limiting slow motion along
an edge attached to the vertex is, respectively, directed from or to
the vertex).

Note that, since we assume that $\lim\li_{|x|\ra
\infty}F(x)=\infty$, at least one local maximum of $G(x)$ is
available on each connected component of every level set of $F(x)$.

\

\section{Positive genus level set components}

Consider a slightly more general equation

$$\dot{\widetilde{X}}_t=\grad F (\widetilde{X}_t)\times \mathbf{d}(\widetilde{X}_t)
\ , \ \widetilde{X}_0=x(z) \ , \eqno(5.1)$$ where the initial point
$x(z)$ belongs to one of the connected components $M=M(z)$ of the
level set $\{x\in \R^3: F(x)=z\}$. As before, we assume that $F(x)$
is smooth enough, $\lim\li_{|x|\ra \infty}F(x)=\infty$, and $\grad
F(x)\neq 0$ for $x\in M$, so that $M$ is a compact connected
orientable two-dimensional surface in $\R^3$.

The vector field $\mathbf{d}(x)$, $x\in \R^3$, is assumed to be
smooth and the vector field $\grad F (x)\times \mathbf{d}(x)$ has,
at most, a finite number of rest points on $M$. Moreover, assume
that
$$\grad \times \mathbf{d}(x)=0 \text{ for } x\in M \ .$$ Note that
in the case of equation (2.1), $\mathbf{d}(x)=\grad G(x)$, and the
last assumption is satisfied.

We will make use of the following

\

\textbf{Lemma 5.1.} \textit{The measure on $M(z)$ with the density
with respect to the surface area proportional to $\dfrac{1}{|\grad
F(x)|}$ is invariant for the flow (5.1) on $M(z)$.}

\

\textit{Proof}. Let us consider an auxiliary system
$$\dot{\widetilde{\widetilde{X}}}_t=\dfrac{\grad F (\widetilde{\widetilde{X}}_t)}
{|\grad F (\widetilde{\widetilde{X}}_t)|}\times
\mathbf{d}(\widetilde{\widetilde{X}}_t) \ , \
\widetilde{\widetilde{X}}_0=x(z) \ ,$$ which is a time change of
system (5.1). Take any closed non self-intersecting curve $\gm$ on
$M$ bounding a region $D(\gm)$ on $M$. Let the unit vector field
$\mathbf{e}_1$ be outward normal to $\gm$, but tangent to $M$. Let
$\mathbf{e}_3=\dfrac{\grad F}{|\grad F|}$. Let the unit vector field
$\mathbf{e}_2$ be tangent to $\gm$ and $M$:
$\mathbf{e}_2=\mathbf{e}_3\times \mathbf{e}_1$. We have

$$\begin{array}{l}
\play{\oint_{\gm}\left(\dfrac{\grad F}{|\grad F|}\times
\mathbf{d}\right) \cdot \mathbf{e}_1 dl}
\\
\play{=-\oint_{\gm}\left(\mathbf{d}\times \mathbf{e}_3\right) \cdot
\mathbf{e}_1 dl}
\\
\play{=-\oint_{\gm}\mathbf{d}\cdot \left(\mathbf{e}_3 \times
\mathbf{e}_1 \right)dl =-\oint_{\gm}\mathbf{d}\cdot \mathbf{e}_2dl
=-\iint_{D(\gm)}\grad \times \mathbf{d}dm=0 \ .}\end{array}$$ (Here
$dm$ is the area element on $M$.)

Therefore the flow of the auxiliary system
$\widetilde{\widetilde{X}}_t$ is incompressible (divergence-free) on
$M$. Thus the standard surface area (induced by the metric element
in $\R^3$) is invariant for $\widetilde{\widetilde{X}}_t$. Since
$\widetilde{X}_t$ is a time change of $\widetilde{\widetilde{X}}_t$
with a factor $|\grad F(x)|$, we see that the measure on $M(z)$ with
the density proportional to $\dfrac{1}{|\grad F(x)|}$ is invariant
for flow (5.1) on $M(z)$. $\square$

\

The topological structure of a compact two-dimensional orientable
connected manifold $M$ is uniquely determined by its genus. If the
genus of $M$ is zero, the condition $\grad \times \mathbf{d}(x)=0$
for $M$ implies that $\mathbf{d}(x)=\grad G(x)$ for a smooth
function $G(x)$. Perturbation theory for such systems was considered
in Sections 3 and 4.

But in the case when $M$ has higher genus, situation is more
complicated. Let us consider, for example, the case when the genus
of $M$ is $1$ so that $M=\torus$ is a two-dimensional torus. The
general structure of an area preserving flow on a torus is described
in [2] (Also see [19, Theorem 3.1.7]. Here not exactly the area is
preserved, but a measure with strictly positive and bounded density.
Then the structure of the trajectories is similar to the case of
area-preserving systems on $M$): There exist finitely many domains
$U_k\subset \torus$ ($k=1,...,n$), bounded by the separatrices of
the flow, such that the trajectories of the dynamical system (5.1)
in each $U_k$ behaves as in a part of the plane: they are either
periodic or tend to a point where the vector field is equal to zero.
Outside of the domains $U_k$ the trajectories form one ergodic
class. Let this ergodic class be
$\cE=\torus\setminus(\bigcup\li_{k=1}^n \overline{U}_k)$ (here and
below $\overline{U}_k$ is the closure of $U_k$). Within each $U_k$
the system (5.1) behaves like a standard Hamiltonian system with a
Hamiltonian $H_k$. For brevity let us assume that each $U_k$
contains only one maxima or minima of $H_k$ and no saddles (the case
when there is a saddle can be resolved using the results of previous
sections). Let us denote the maxima or minima of $H_k$ in $U_k$ by
$M_k$. Let $A_k$ be the saddles of $(5.1)$ on $\torus$: $A_k$ is
situated on the boundary of $U_k$. Let us introduce a family of
functions $h_k(x)=H_k(x)-H_k(A_k)$ when $x\in U_k$ and $h_k(x)=0$
when $x\in \torus\setminus\overline{U}_k$, $k=1,...,n$. Let the set
$\{x\in \overline{U}_k; h_k(x)=h_k\}$ be $\gm_k(h_k)$. We notice
that $\gm_k(0)$ is the separatrix bounding $U_k$ and containing
$A_k$.

Identify all points of the ergodic class $\cE$ as well as the points
belonging to each level set of each function $H_k(x)$, $x\in
\overline{U}_k$. Let $\iiY$ be the identification mapping. Then
$\iiY(M)$, in the natural topology, is homeomorphic to a graph
$\graph$. This graph is a tree, and $\iiY$ maps the entire ergodic
class $\cE$ to the root of the graph which is denoted by $O$. Let
$\gm_k(h)=\{x\in \overline{U}_k: h_k(x)=h\}$. Define a metric
$\rho(y_1,y_2)$ on $\graph$ as follows: If $y_1=\iiY(\gm_k(h_1))$,
$y_2=\iiY(\gm_l(h_2))$, put $\rho(y_1,y_2)=|h_1-h_2|$ for $k=l$, and
$\rho(y_1,y_2)=\rho(y_1,O)+\rho(O,y_2)$ if $k \neq l$. In this way
the region $\overline{U}_k$ will be mapped into a segment $I_k$ of
the form either $[0, h_k(M_k)]$ (if $M_k$ is a maximum) or
$[h_k(M_k),0]$ (if $M_k$ is a minimum). All these segments $I_k$
serve as edges of our graph $\graph$ and they share the common root
$O$. Every point $y=\iiY(x)$ on $\graph \setminus O$ can be given a
coordinate $(k, h_k)$ where $k$ is the number of the edge containing
$y$ and $h_k=h_k(x)$. In this way our mapping $\iiY$ is explicitly
written as $\iiY(x)=O$ if $x\in \overline{\cE}$ and
$\iiY(x)=(k,h_k(x))$ if $x\in U_k$.

Let us now introduce a deterministic perturbation and a stochastic
regularization to our system (5.1). After the time change $t \mapsto
\dfrac{t}{\ve}$, our perturbed system has the form

$$\dot{X}_t^{\ve,\dt}=\dfrac{1}{\ve}\grad F(X_t^{\ve,\dt}) \times \mathbf{d }(X_t^{\ve,\dt})+
\grad F (X_t^{\ve,\dt})\times\mathbf{ p} (X_t^{\ve,\dt})+\dt
\sm(X_t^{\ve,\dt})\circ \dot{W}_t \ , \ X_0^{\ve,\dt}=x_0(z) \ .
\eqno(5.2)$$

Here $\mathbf{p}(\bullet)$ is a smooth vector field in $\R^3$ and
$\sm$ is the same matrix defined in Section 4. We remind the reader
that $\sm^T \grad F =0$ and $a=(a_{ij})=\sm\sm^T$ is the diffusion
matrix. We also recall that we have the non-degeneracy conditions of
$a$ on $M$: $\mathbf{e}\cdot(a(x)\mathbf{e})\geq
\underline{a}|\mathbf{e}|^2$ for some $\underline{a}>0$ and all
$\mathbf{e}$ such that $\mathbf{e}\cdot \grad F=0$. The process
$X_t^{\ve,\dt}$ lives on the surface $M$.

Let us define a strong Markov process $Y_t^\dt$ on $\graph$ as the
diffusion process on $\graph$ governed by a generator $A$ such that,
at each interior point $(k,h_k)$ of an edge $I_k$,
$Af(k,h_k)=\overline{L}_kf(k,h_k)$, where

$$\overline{L}_k f(k, h_k)=\dfrac{1}{T_k(h_k)}\left(a_k(h_k)+\dfrac{\dt^2}{2}a_{1,k}(h_k)
+\dfrac{\dt^2}{2}a_{2,k}(h_k)\right)\dfrac{\pt f}{\pt
h_k}+\dfrac{\dt^2}{2}\dfrac{1}{T_k(h_k)}b_k(h_k)\dfrac{\pt^2 f}{\pt
h_k^2} \ , \eqno(5.3)$$

with

$$\begin{array}{l}
\play{a_k(h_k)=\oint_{\gm_k(h_k)}\grad H_k\cdot (\grad F \times
\mathbf{p})\dfrac{dl}{|\grad F \times \mathbf{d}|} \ ,}
\\
\play{a_{1,k}(h_k)=\oint_{\gm_k(h_k)}\grad H_k \cdot \mathbf{\Sm}
\dfrac{dl}{|\grad F \times \mathbf{d}|} \ ,}
\\
\play{a_{2,k}(h_k)=\oint_{\gm_k(h_k)}\sum\li_{i,j=1}^3a_{ij}\dfrac{\pt^2
H_k}{\pt x_i\pt x_j}\dfrac{dl}{|\grad F \times \mathbf{d}|} \ ,}
\\
\play{b_k(h_k)=\oint_{\gm_k(h_k)}|(\grad
H_k)^T\sm|^2\dfrac{dl}{|\grad F \times \mathbf{d}|} \ ,}
\end{array} \eqno(5.4)$$
and $$T_k(h_k)=\oint_{\gm_k(h_k)}\dfrac{dl}{|\grad F \times
\mathbf{d}|}$$ is the period of one rotation along $\gm_k(h_k)$.
Here the vector $\mathbf{\Sm}$ is the same vector as in Section 4.

The domain $D(A)$ of $A$ consists of those functions $f$ that are
continuous on $\graph$ and have the following properties.

$\bullet$ Function $f$ is twice continuously differentiable in the
interior of each of the edges.

$\bullet$ We have the one sided limits $\lim\li_{h_k \ra
0}\overline{L}_kf(k,h_k)$ and $\lim\li_{h_k\ra
h_k(M_k)}\overline{L}_kf(k,h_k)$ at the endpoints of each of the
edges. The values of the limit $q=\lim\li_{h_k\ra
0}\overline{L}_kf(k,h_k)$ are the same for all the edges.

$\bullet$ The following gluing condition is satisfied at $O$:

$$\sum\li_{k=1}^n(\pm)\bt_k\lim\li_{h_k\ra 0}\dfrac{\pt f}{\pt h_k}(k,h_k)
=q \ , \eqno(5.5)$$ with sign $+$ if $A_k$ is a local minimum of
$H_k$ restricted on $U_k$ and sign $-$ otherwise. Here
$$\bt_k=\dfrac{1}{\lb(\cE)}\oint_{\gm_k(0)}|(\grad H_k)^T \sm|^2 \dfrac{dl}{|\grad F \times \mathbf{d}|}$$
with $$\lb(\cE)=\iint_{\cE}\dfrac{dm}{|\grad F|}\ . $$ $\ (\text{
Here } dm \text{ is the area element on } M .)$ These conditions
define the process $Y_t^\dt$ on $G$ in a unique way.

We have the following

\

\textbf{Theorem 5.1.} \textit{The process}
$Y_t^{\ve,\dt}=\iiY(X_t^{\ve,\dt})$ \textit{converges weakly in the
space of continuous trajectories $[0,T]\ra \graph$ as $\ve \da 0$ to
$Y_t^\dt$.}

\

The \textit{proof} of this theorem is an application of Theorem 1 of
[6]. To be precise, in formula (5) of [6], we set $\kp=\dt^2$ ,
$v(X_t^{\ve,\dt})=\grad F(X_t^{\ve,\dt})\times \mathbf{d}
(X_t^{\ve,\dt})$, $\bt(X_t^{\ve,\dt})=\grad F(X_t^{\ve,\dt})\times
\mathbf{p}(X_t^{\ve,\dt})$,
$u(X_t^{\ve,\dt})=\dfrac{\mathbf{\widetilde{c}}(X_t^{\ve,\dt})}{2}$
(a term which comes from the Stratonovich integral),
$\sm(X_t^{\ve,\dt})=\sm(X_t^{\ve,\dt})$. Furthermore, we can write
down the generator $L$ of $X_t^{\ve,\dt}$ in self-adjoint form
$$Lu=\dfrac{\dt^2}{2}\sum\li_{i=1}^3\dfrac{\pt}{\pt x_i}
\left(\sum\li_{j=1}^3a_{ij}\dfrac{\pt u}{\pt x_j}\right)
+\left(\dfrac{1}{\ve}\grad F \times \mathbf{d}+\grad F \times
\mathbf{p}-\dfrac{\dt^2}{2}\mathbf{\Pi}\right)\cdot \grad u \ .$$
Here $\mathbf{\Pi}$ is a $3$-vector with the $i$-th component
$\Pi_i=\sum\li_{j,k=1}^3\dfrac{\pt \sm_{kj}}{\pt x_k}\sm_{ij}$.
Notice that since $\sm^T\grad F=0$, we have $\grad F \cdot
\mathbf{\Pi}=\play{\sum\li_{j.k=1}^3}\dfrac{\pt \sm_{kj}}{\pt
x_k}\sum\li_{i=1}^3\sm_{ij}\dfrac{\pt F}{\pt x_i}=0$. Also notice
that since we have checked the fact that $F(X_t^{\ve,\dt})$ is a
constant of motion, It\^{o}'s formula imply $LF(x)=0$. Therefore, we
have $\mathcal{L}F(x)=0$ where
$$\mathcal{L}u=\dfrac{\dt^2}{2}\sum\li_{i=1}^3\dfrac{\pt}{\pt x_i}
\left(\sum\li_{j=1}^3a_{ij}\dfrac{\pt u}{\pt x_j}\right) \ .
$$
From here we see that the auxiliary process $\mathcal{X}_t^{\ve,\dt}
\ , \ \mathcal{X}_0^{\ve,\dt}=x_0(z)$ corresponding to the operator
$\mathcal{L}$ lives on the surface $M=M(z)$. Since $\mathcal{L}$ is
self-adjoint in $\R^3$, the (degenerate) process
$\mathcal{X}_t^{\ve,\dt}$ has an invariant measure proportional to
$\R^3$ Lebesgure measure. This implies, that the process
$\mathcal{X}_t^{\ve,\dt}$, viewed as a non-degenerate diffusion
process on $M$, has a unique invariant measure with density
proportional to $\dfrac{1}{|\grad F(x)|}$ (with respect to the
surface area element $dm$ on $M$). Since we have checked that the
deterministic flow (5.1) on $M$ also has an invariant measure with
density proportional to $\dfrac{1}{|\grad F(x)|}$, we see that the
auxiliary process $\mathfrak{X}_t^{\ve,\dt}$,
$\mathfrak{X}_t^{\ve,\dt}=x_0(z)$ governed by the operator
$$\mathfrak{L}u=\dfrac{\dt^2}{2}\sum\li_{i=1}^3\dfrac{\pt}{\pt x_i}
\left(\sum\li_{j=1}^3a_{ij}\dfrac{\pt u}{\pt x_j}\right)
+\left(\dfrac{1}{\ve}\grad F \times \mathbf{d}\right)\cdot \grad u
$$ is a non-degenerate diffusion process on $M$ with a unique
invariant measure which has a density proportional to
$\dfrac{1}{|\grad F(x)|}$. This fact, together with the standard
method of absolutely continuous change of measure (see [11] and
compare with Appendix.2), allow us to calculate the gluing condition
(5.5).

\

Since the small random perturbation term $\dt\sm\circ \dot{W}_t$ in
(5.2) is only introduced as a regularization, we must study the
limit of $Y_t^\dt$ as $\dt \da 0$. It follows from the same argument
as in Section 3 of [6] that the limiting process $Y_t$ should be
described as follows. Let

$$\overline{\psi}_k=2\oint_{\gm_k(0)}\grad H_k \cdot (\grad F
\times \mathbf{p})\dfrac{dl}{|\grad F \times \mathbf{d}|}\neq 0 \
.$$

Let $s_k$, $1\leq k \leq n$, take values $0$ and $1$. We set $s_k=1$
if $\overline{\psi}_k>0$ and $M_k$ is a local maximum of $H_k$ as
well as if $\overline{\psi}_k<0$ and $M_k$ is a local minimum of
$H_k$. Otherwise we set $s_k=0$. Let
$$r_k=\dfrac{s_k|\overline{\psi}_k|}{2\lb(\cE)} \ ,  1\leq k \leq n \ .$$
Then we can describe $Y_t$ as follows.

$\bullet$ The process $Y_t$ is a strong Markov process with
continuous trajectories.

$\bullet$ If $Y_0=O$, where $O$ is the root of $\graph$, then the
process spends a random time $\tau$ in $O$. There is a random
variable $\xi$ that is independent of $\tau$, taking values in the
set $\{1,...,n\}$, such that $Y_t\in I_{\xi}$ for $t>\tau$. If
$s_k=0$ for all $k$, $1\leq k \leq n$ then $\tau=\infty$. If $s_k=1$
for some $k$ then $\tau$ is distributed as an exponential random
variable with expectation $\sum\li_{k=1}^n r_k$. If $s_k=1$ for some
$k$ then

$$\Prob(Y_t\in I_k, t>\tau)=\dfrac{r_k}{\sum\li_{k=1}^n r_k} \ .$$

$\bullet$ If $Y_0\in \text{Int} I_k$ then
$$\dfrac{dY_t}{dt}=\overline{B}_k(Y_t)$$ for $t<\sm$ where $\sm=\inf(t:
Y_t=0)$ and
$\overline{B}_k(h_k)=\dfrac{\overline{\psi}_k(h_k)}{2T_k(h_k)}$.

\

\textbf{Theorem 5.2.} \textit{As $\dt\da 0$, the process $Y_t^\dt$
converges weakly in the space of continuous trajectories $[0,T]\ra
\graph$, to the process $Y_t$.}

\

The \textit{proof} is an application of Theorem 2 in [6]. (See the
explanation in the proof of Theorem 5.1.)

\

In the more general situation when the surface $M$ has higher genus,
the situation is similar (compare with [7]). In particular,
corresponding graph may be not a tree; it can have more than one
special vertices where the limiting Markov process spends random
time with exponential distribution; transitions between those
special vertices are possible.

\

\section{Multiscale perturbations}

Equation (2.1) has two first integrals $F(x)$ and $G(x)$. These
integrals may have different nature and their perturbations may have
different order. Consider the case when the perturbed system has the
form

$$\begin{array}{l}
\dot{X}_t^{\ve,\kp}=\grad F (X_t^{\ve,\kp})\times \grad
G(X_t^{\ve,\kp})+
\sqrt{\kp}\sm_1(X_t^{\ve,\dt})*\dot{W}_t^1+\sqrt{\ve}\sm_2(X_t^{\ve,\dt})*\dot{W}_t^2
\ ,  \\  \ve,\kp>0 \ , \  X_0^{\ve,\kp}=x\in M\subset\{y\in \R^3:
F(y)=z\} \ ,
\end{array} \eqno(6.1)
$$
where $M$ is a connected component of the level set $\{F(x)=z\}$;
$\sm_1(x)$ and $\sm_2(x)$ are $3\times 3$-matrices; $\dot{W}_t^1$
and $\dot{W}_t^2$ are independent white noises in $\R^3$. Put
$a_1(x)=\sm_1(x)\sm_1^T(x)$, $a_2(x)=\sm_2(x)\sm_2^T(x)$. Sign "$*$"
in the stochastic terms means that the stochastic integrals are
defined in such a way, that the generator of the process
$X_t^{\ve,\kp}$ is as follows

$$L^{\ve,\kp}u(x)=(\grad F (x) \times \grad G(x)) \cdot \grad u(x) +
\dfrac{\kp}{2}\text{div}(a_1(x)\grad
u(x))+\dfrac{\ve}{2}\text{div}(a_2(x)\grad u(x))\ . \ \eqno(6.2)$$

We assume that $a_1(x)\grad F(x)=0$ and $\mathbf{e}\cdot
(a_1(x)\mathbf{e})\geq \underline{a}_1 |\mathbf{e}|^2$ for each
$\mathbf{e}$ such that $\mathbf{e}\cdot \grad F(x)=0$,
$\underline{a}_1$ is a positive constant. The matrix $a_2(x)$ is
assumed to be non-degenerate. The assumptions concerning $a_1(x)$
imply that the process $X_t^{0,\kp}$ moves on the surface $M$:
$\Prob\{X_t^{0,\kp}\in M\}=1$. This follows directly from the
It\^{o} formula (we refer the reader to the proof of Theorem 5.1 in
Section 5, where we did a similar calculation). Moreover, the
process $X_t^{0,\kp}$ on $M$ is non-degenerate. This implies that,
for any $\kp>0$, the process $X_t^{0,\kp}$ has on the compact
manifold $M$ (we assume that $\lim\li_{|x|\ra \infty}F(x)=\infty$) a
unique invariant measure. On the other hand, the drift in (6.2) is
divergence-free and the main part is formally self-adjoint.
Therefore the Lebesgue measure is invariant for the process
$X_t^{\ve,\kp}$, and in particular for $X_t^{0,\kp}$, in $\R^3$.
This implies that $\dfrac{C}{|\grad F(x)|}$,
$\play{C=\left(\int_M\dfrac{dm}{|\grad F(x)|}\right)^{-1}}$, where
$dm$ is the surface area on $M$, is the density of the unique
invariant measure of $X_t^{0,\kp}$ on $M$ for each $\kp>0$.

Assume that $0<\ve<<\kp<1$. This means that we have relatively large
perturbations of the first integral $G(x)$ and much smaller
perturbations of $F(x)$. On the time intervals of order
$\dfrac{1}{\kp}$, one can omit the term
$\sqrt{\ve}\sm_2(x)*\dot{W}_t^2$ in (6.1): the first integral
$F(X_t^{\ve,\kp})$ does not change on such intervals as $0\leq
\ve<<\kp<<1$, and the evolution of $G(X_t^{\ve,\kp})$ asymptotically
coincides with the evolution of $G(X_t^{0,\kp})$ and can be
described using the results of Section 4.

But on time intervals of order $\dfrac{1}{\ve}>\dfrac{1}{\kp}$, the
situation is different. Consider process
$\widehat{X}_t^{\ve,\kp}=X_{t/\ve}^{\ve,\kp}$. The process
$\widehat{X}_t^{\ve,\kp}$ is governed by the generator
$\dfrac{1}{\ve}L^{\ve,\kp}=\widehat{L}^{\ve,\kp}$. It has a fast and
a slow components as $\ve \da 0$. The fast component of the process
$\widehat{X}_t^{\ve,\kp}$ can be approximated by the process
$\widehat{\widehat{X}}_t^{\ve,\kp}$ corresponding to the generator
$$\widehat{\widehat{L}}^{\ve,\kp}u(x)=\dfrac{1}{\ve}(\grad F(x)\times \grad G(x)) \cdot
\grad u + \dfrac{\kp}{2\ve}\text{div}(a_1(x)\grad u) \ . $$ The
process $\widehat{\widehat{X}}_t^{\ve,\kp}$ lives on the surface $M$
and, up to a simple time change $t \ra \dfrac{t}{\ve}$, coincides
with $X_t^{0,\kp}$. In particular, it has the same invariant density
$C|\grad F(x)|^{-1}$.

To describe the slow component of $\widehat{X}_t^{\ve,\kp}$, one
should introduce a graph. Identify points of each connected
component of every level set of the function $F(x)$. Let $\iY$ be
the identification mapping. Then the set $\iY(\R^3)$ is homeomorphic
to a graph provided with the natural topology which we denote by
$\Gm$.

Note that all connected components of level sets not containing
critical points of $F(x)$ are two-dimensional compact (we assume
that $\lim\li_{|x|\ra \infty}F(x)=\infty$ manifolds). Each local
maximum or minimum of $F(x)$ corresponds to an exterior vertex
belonging just to one edge. The saddle points correspond to the
interior vertices. Unlike in the case of generic functions of two
variables, not every interior vertex belongs to three edges: If $O$
is a saddle point of $F(x)$, the surface $\{y\in \R^3: F(y)=F(O)\}$
divides each small neighborhood of $O$ in three parts. But two of
these parts, in the case of functions of three variables can come
together far from $O$ (compare with [12]). One can introduce a
global coordinate system on $\Gm$: Number the edges of $\Gm$. Then
each point $y\in \Gm$ can be identified by two numbers $k$ and $z$,
where $k$ is the number of an edge containing $y$ and
$z=F(\iY^{-1}(y))$.

The slow component of $X_t^{\ve,\kp}$ is the (not Markovian, in
general) process $\iY(\widehat{X}_t^{\ve,\kp})=Y_t^{\ve,\kp}$ on
$\Gm$.

Define a diffusion process $Y_t$ on $\Gm$ which inside each edge
$I_k\subset \Gm$ is governed by an ordinary differential operator
$\overline{L}_k=\dfrac{1}{2T_k(z)}\dfrac{d}{dz}(\overline{a}_k(z)\dfrac{d}{dz})$,
where

$$T_k(z)=\int_{\iY^{-1}(k,z)}\dfrac{dm}{|\grad F(x)|} \ ,
\ \overline{a}_k(z)=\int_{G(k,z)}\text{div}(a_2\grad F(x))dx \ ,
\eqno(6.3)$$ where $G(k,z)\subset \R^3$ is the domain bounded by the
surface $\iY^{-1}(k,z)$; $a_2(x)=\sm_2(x)\sm_2^{T}(x)$, $dm$ is the
area element on $\iiY^{-1}(k,z)$.

The operators $\overline{L}_k$ define the process $Y_t$ inside the
edges. To define the behavior of $Y_t$ at the vertices, we describe
the domain $D_A$ of the generator of $Y_t$ (see Ch.8 in [14]). We
say that a continuous on $\Gm$ and smooth inside the edges function
$f\in D_A$ if and only if the following holds.

$\bullet$ The function defined inside the edges by the formula
$\overline{L}_kf(k,z)$ can be extended to a continuous on the whole
graph function.

$\bullet$ If edges $I_{i_1}$, $I_{i_2}$, $I_{i_3}$ are attached to
an interior vertex $O$, then $$\sum\li_{k=1}^3
(\pm)\overline{a}_{i_k}(O)D_kf(O)=0 \ , $$ where
$\overline{a}_{i_k}(O)=\lim\li_{z\ra F(O)}\overline{a}_{i_k}(z)$
($\overline{a}_i(z)$ is defined by (6.3)), and
$$D_kf(O)=\lim\li_{z\ra F(O)}\dfrac{f(k,z)-f(k,F(O))}{z-F(O)}$$
(compare with [12]). The sign convention in the gluing condition is
as follows: Let $\iiY^{-1}(I_{i_1})$ belong to the set $\{x\in \R^3:
F(x)\geq F(O)\}$, and $\iiY^{-1}(I_{i_2})$,
$\iiY^{-1}(I_{i_3})\subset \{x\in \R^3: F(x)\leq F(O)\}$. Then sign
$+$ should be taken in front of $\overline{a}_{i_1}(O)$ and sign $-$
in front of $\overline{a}_{i_2}(O)$ and $\overline{a}_{i_3}(O)$.

$\bullet$ If just two edges $I_{i_1}$ and $I_{i_2}$ are attached to
an interior vertex $O$, then $D_{i_1}f(O)=D_{i_2}f(O)$.

\

For functions $f(k,z)$ with these properties,
$Af(k,z)=\overline{L}_k f(k,z)$. These conditions define the Markov
process $Y_t$ on $\Gm$ in a unique way. Exterior vertices are
inaccessible for $Y_t$.

\

\textbf{Theorem 6.1.} \textit{The process
$Y_t^{\ve,\kp}=\iiY(\widehat{X}_t^{\ve,\kp})$ converges weakly in
the space of continuous functions $[0,T]\ra \Gm$ for each finite
$T>0$ as $\ve \da 0$ to the (independent of $\kp$ and $\sm_1(x)$)
process $Y_t$ defined above.}

\

The \textit{proof} of this statement follows from Theorem 2.1 of
[12]. We omit the details. Using the absolute continuity arguments
which we mentioned earlier one can consider more general
perturbations in (6.1).

\

\section*{Appendix}

1. We provide here the proof of Lemma 3.1. By a similar calculation
as we did before stating Lemma 3.1 we have

$$\left|\dfrac{L(a(z),b(z))}{L(b(z),c(z))}-
\dfrac{\text{Area}(\square_1)}{\text{Area}(\square_2)}\right|
<\dt_1(\lb)+\dt_2(\ve) \ . \eqno(\text{A}.1.1)$$ (Here and below we
use symbol $\dt_k(\mu)$ to denote a positive quantity which goes to
zero as the parameter $\mu \da 0$.)

We can also check, by mean value theorem, that

$$\left|\dfrac{\text{Area}(\square_1)}{\text{Area}(\square_2)}
-\dfrac{\play{\iint_{\square_1}\left|\dfrac{1}{\ve}\grad F \times
\grad G +\grad F \times (\grad F \times
\mathbf{b})\right|dm}}{\play{\iint_{\square_2}\left|\dfrac{1}{\ve}\grad
F \times \grad G +\grad F \times (\grad F \times
\mathbf{b})\right|dm}}\right|<\dt_3(\lb\ve) \ .
\eqno(\text{A}.1.2)$$

By (3.3), it is easy to check that

$$\left|
\dfrac{\play{\iint_{\square_1}\left|\dfrac{1}{\ve}\grad F \times
\grad G + \grad F \times (\grad F \times \mathbf{b})\right|dm}}
{\play{\iint_{\square_2}\left|\dfrac{1}{\ve}\grad F \times \grad G +
\grad F \times (\grad F \times
\mathbf{b})\right|dm}}-\dfrac{\play{\iint_{\cS_1(z)}\grad \times
(\grad F \times \mathbf{b})\cdot \mathbf{n}dm}}
{\play{\iint_{\cS_2(z)}\grad \times (\grad F \times \mathbf{b})\cdot
\mathbf{n}dm}}\right|<\dt_4(\lb) \ . \eqno(\text{A}.1.3)$$

By using the averaging principle, it is possible to show that the
ratio $\dfrac{\text{Area}(\square_1)}{\text{Area}(\square_2)}$ is
asymptotically preserved along the flow of (2.7) (compare with [4]).
Therefore we can take $\square_1$ and $\square_2$ as close to the
separatrices hitting and exiting $O_2'(z)$ as we wish. This fact,
together with the estimates (A.1.1)-(A.1.3), imply our Lemma 3.1, by
letting first $\lb \da 0$ and then $\ve \da 0$.

\

2. We explain here the missing details in the proof of Theorem 4.1.
As we have explained in that proof, our process $Z_t^{\ve,\dt}$
satisfies the equation

$$\dot{Z}_t^{\ve,\dt}=\dfrac{1}{\ve}\kp(Z_t^{\ve,\dt})\overline{\grad}\cG(Z_t^{\ve,\dt})
+\vec{\bt}_1(Z_t^{\ve,\dt})+
\dfrac{\dt^2}{2}\mathbf{\widetilde{c}}(Z_t^{\ve,\dt})+\dt\widetilde{\sm}(Z_t^{\ve,\dt})
\dot{\widetilde{W}}_t \ , \
Z_0^{\ve,\dt}=z_0=(f_1(x_0(z)),f_2(x_0(z))) \ .
\eqno(\text{A}.2.1)$$

Here the term
$\dfrac{\dt^2}{2}\mathbf{\widetilde{c}}(Z_t^{\ve,\dt})$ comes from
the Stratonovich integral in (4.9).

As before, we can identify the connected components of the level
sets of the Hamiltonian $\cG$ to obtain a graph $\Gm$. Let $\iiY$ be
the identification mapping. Let us use the same symbols to denote
vertices and edges as those we use for the graph corresponding to
$X_t^{\ve,\dt}$ (see Section 2).

System (A.2.1), by a non-singular time change, can be reduced to a
perturbed Hamiltonian system with Hamiltonian $\cG$. The form of the
operators governing the limiting diffusion inside the edges is
obtained by standard averaging. To get the gluing conditions, we
first consider an auxiliary process

$$\dot{\widehat{Z}}_t^{\ve,\dt}=
\dfrac{1}{\ve}\kp(\widehat{Z}_t^{\ve,\dt})
\overline{\grad}\cG(\widehat{Z}_t^{\ve,\dt}) + \dt
\widetilde{\sm}(\widehat{Z}_t^{\ve,\dt})\dot{\widetilde{W}}_t \ , \
\widehat{Z}_0^{\ve,\dt}=z_0 \ . \eqno(\text{A}.2.2)$$

Such a process, by a non-singular time change, is equivalent to a
perturbed Hamiltonian system which has Lebesgue measure as its
invariant measure. Using this fact, via a standard proof of [14,
Chapter 8, Section 6], we conclude that the gluing condition for the
weak limit of $\iiY(\widehat{Z}_t^{\ve,\dt})$ as $\ve \da 0$ at
vertex $O_2(z)$ is given by the coefficients

$$\bt_{2,i}=\oint_{f(C_i(0,z))}|(\grad \cG)^T\widetilde{\sm}|^2\dfrac{dl_z}{|\vec{\bt}|} \ ,$$
for $i=1,2,3$. Here $\vec{\bt}=\kp\overline{\grad} \cG$.

The measure $\widehat{\mu}^{\ve,\dt}$ corresponding to
$\widehat{Z}_t^{\ve,\dt}$ ($0\leq t \leq T$) is related to the
measure $\mu^{\ve,\dt}$ corresponding to $Z_t^{\ve,\dt}$ ($0\leq t
\leq T$) via the Girsanov formula

$$\begin{array}{l}
\play{\dfrac{d\mu^{\ve,\dt}}{d\widehat{\mu}^{\ve,\dt}}=I_{0T}^{\ve,\dt}=
\exp\left\{\dfrac{1}{\dt}\int_0^T\sm^{-1}
(\widehat{Z}_t^{\ve,\dt})[\vec{\bt}_1(\widehat{Z}_t^{\ve,\dt})+
\dfrac{\dt}{2}\mathbf{\widetilde{c}}(\widehat{Z}_t^{\ve,\dt})]\cdot
d\widetilde{W}_t - \right.}
\\
\play{\left. \ \ \ \ \ \ \ \ \ \ \ \ \ \ \ \ \ \ \ \ \ \ \ \ \ \ \ -
\dfrac{1}{2\dt^2}\int_0^T
\left|\sm^{-1}(\widehat{Z}_t^{\ve,\dt})[\vec{\bt}_1(\widehat{Z}_t^{\ve,\dt})+\dfrac{\dt}{2}
\mathbf{\widetilde{c}}(\widehat{Z}_t^{\ve,\dt})]\right|^2dt
\right\}} \ .
\end{array}$$

\

\textbf{Lemma A.2.1.} \textit{There exist constants $A_1>0$, $T_0>0$
such that $\E_{z_0}(I_{0T}^{\ve,\dt}-1)^2\leq A_1 T$ for all
$T<T_0$.}

\

To prove this lemma, we first apply It\^{o}'s formula to
$(I_{0T}^{\ve,\dt}-1)^2$ and taking expected value. After that we
use the fact that

$$\begin{array}{l}
\E_{z_0}\play{\exp\left\{\dfrac{2}{\dt}\int_0^T\sm^{-1}
(\widehat{Z}_t^{\ve,\dt})[\vec{\bt}_1(\widehat{Z}_t^{\ve,\dt})+
\dfrac{\dt^2}{2}\mathbf{\widetilde{c}}(\widehat{Z}_t^{\ve,\dt})]\cdot
d\widetilde{W}_t - \right.}
\\
\play{\left. \ \ \ \ \ \ \ \ \ \ \ \ \ \ \ \ \ \ \ \ \ \ \ \ \ \ \ -
\dfrac{2}{\dt^2}\int_0^T
\left|\sm^{-1}(\widehat{Z}_t^{\ve,\dt})[\vec{\bt}_1(\widehat{Z}_t^{\ve,\dt})+\dfrac{\dt^2}{2}
\mathbf{\widetilde{c}}(\widehat{Z}_t^{\ve,\dt})]\right|^2dt
\right\}=1}
\end{array}$$
and the Cauchy-Schwarz inequality. The proof is essentially the same
as that of Lemma 2.3 in [11].

\

For small $\lb>0$ we let

$$D_2(\lb)=\{x\in \R^2: \cG(x)\in [-\lb,\lb]\} \ .$$

For $i=1,3$, we let

$$D_i(\lb)=\{x\in \R^2: \cG(f(O_i(z)))\leq \cG(x)\leq
\cG(f(O_i(z)))+\lb \ , \ x \text{ is in the well } D_i(0,z) \text{
containing } O_i(z)\} .$$

For $k=1,2,3$ we let

$$\tau_k^{\ve,\dt}(\lb)=\inf \{t>0, Z_t^{\ve,\dt}\not \in D_k(\lb)\} \ .$$

We have

\

\textbf{Lemma A.2.2.} \textit{For any positive $\mu>0$ and $\kp>0$
there exists $\lb_0>0$ such that for $0\leq \lb < \lb_0$ for
sufficiently small $\ve$ and all $x\in D_2(\lb)$}

$$\E_{z_0}\int_0^{\tau_2^{\ve,\dt}(\lb)}\exp(-\mu t)dt<\kp \lb \ ,$$
\textit{and for all $x\in D_i(\lb)$ ($i=1,3$) we have }

$$\E_{z_0}\int_0^{\tau_i^{\ve,\dt}(\lb)}\exp(-\mu t)dt<\kp \ .$$

\

The \textit{proof} of this Lemma is based on corresponding estimates
for the process $\widehat{Z}_t^{\ve,\dt}$ and Lemma A.2.1. It is
essentially the same as that of Lemma 2.4 in [11].

\

\textbf{Lemma A.2.3.} \textit{Let
$q_i=\dfrac{\bt_{2,i}}{\sum\li_{i=1}^3 \bt_{2,i}}$ where $i=1,2,3$.
We have, for any $\kp>0$ there exist $\lb_0>0$ such that for
$0<\lb<\lb_0$ there exist $\lb'>0$ such that for sufficiently small
$\ve$ we have}
$$\left|\Prob_{z_0}\{Z_{\tau_2^{\ve,\dt}(\lb)}^{\ve,\dt}\in C_{i}((-1)^i\lb,z)\}-q_i\right|<\kp$$
\textit{for all $x\in D_2(\lb')\cup \pt D_2(\lb')$.}

\

The \textit{proof} of this Lemma is also the same as that of Lemma
2.5 in [11].

\

The slow component $\iiY(\widehat{Z}_t^{\ve,\dt})$ of the process
$\widehat{Z}_t^{\ve,\dt}$ converges weakly as $\ve \da 0$ to a
diffusion process $\widehat{Y}_t^\dt$ on $\Gm$. The process
$\widehat{Y}_t^\dt$ is defined by a family of differential
operators, one on each edge of $\Gm$, and by gluing conditions at
the vertices. The operators and gluing conditions were calculated in
Chapter 8 of [14]. The convergence of
$\iiY(\widehat{Z}_t^{\ve,\dt})$ to $\widehat{Y}_t^\dt$ was also
proved in [14].

To find the weak limit of the slow component $\iiY(Z_t^{\ve,\dt})$
of $Z_t^{\ve,\dt}$ as $\ve \da 0$, note that the family
$\iiY(Z_t^{\ve,\dt})$ is weakly compact as $\ve \da 0$. Inside each
edge, the limit is a diffusion process with the generator defined by
the standard averaging principle. The limiting process
$\iiY(Z_t^{\ve,\dt})$ and $\iiY(\widehat{Z}_t^{\ve,\dt})$ inside an
edge, in general, are different. But as it follows from Lemmas
A.2.1-A.2.3, the gluing conditions are the same. This implies that
the family $\iiY(Z_t^{\ve,\dt})$ converges weakly as $\ve \da 0$ and
identifies the limiting process as the process $Y_t^\dt$ in Theorem
4.1.

\

3. We indicate here how to calculate the branching probabilities as
claimed in Theorem 4.2. Let $Y_t^\dt$ be the diffusion process on
graph $\Gm$ described in Theorem 4.1. Let

$$\cE_h(u)=\{v\in \Gm: \rho(u,v)<h\} \text{ for } u\in \Gm ,$$

$$\tau_h^\dt=\min\{t: Y^\dt_t \not \in \cE_h(u)\} \ .$$

\

Let $p_1$ and $p_3$ be defined as in (3.13) and (3.14). We have

\

\textbf{Lemma A.3.1.} \textit{We have, for a small enough $h$, }

$$\lim\li_{\dt \da 0}\Prob_{O_2(z)}(Y^\dt_{\tau^\dt_h}\in I_3)=0 \ ,$$

$$\lim\li_{\dt \da 0}\Prob_{O_2(z)}(Y^\dt_{\tau^\dt_h}\in I_i)=p_i \text{ for } i=1,3 \ .$$

\

To prove this Lemma, we let $u=(g,i)\in \cE_h(O_2(z))$. We set
$v_j^\dt(u)=v_j^\dt(g,i)=\Prob_{(g,i)}\{Y^\dt_{\tau_h^\dt}\in
I_j\}$. The function $v_j^\dt(g,i)$ is the unique continuous
solution of the following problem

$$\left\{\begin{array}{l}
\play{\overline{L}_i v_j^\dt(g,i)=0 \ , \ (g,i)\in
\cE_h(O_2(z))\setminus \{O_2(z)\} \ , \ i=1,2,3 \ ,}
\\
\play{v_j^\dt(g,i)|_{(g,i)\in \pt \cE_h(O_2(z))\cap I_i}=0 \text{
for } i \neq j \ ,}
\\
\play{v_j^\dt(g,j)|_{(g,j)\in \pt \cE_h(O_2(z))\cap I_j}=1 \ ,}
\\
\play{\sum\li_{k=1}^3 (\pm) \bt_{2,k} \lim\li_{g\ra
G(O_2(z))}\dfrac{\pt v_j^\dt}{\pt g}(g,k) = 0 \ .}
\end{array}\right.$$

Here $\overline{L}_i$ are defined in (4.5) and $\bt_{2,i}$ are
defined in (4.7) and (4.8), with $"+"$ sign for $k=2$ and $"-"$ sign
for $k=1,3$. One can solve this problem explicitly and derive the
statement of Lemma A.3.1 similarly to Lemma 2.2 of [4].

\

\textbf{Acknowledgements}: This work is supported in part by NSF
Grants DMS-0803287 and DMS-0854982.

\section*{References}

[1] Arnold V.I., \textit{Mathematical methods of classical
mechanics}, Springer, 1978.

[2] Arnold V.I., Topological and ergodic properties of closed
1-forms with incommensurable periods, \textit{Func. Anal. Appl.}
\textbf{25} (1991), no.2, 81-90.

[3] Athreya A., Freidlin M., Metastability and Stochastic Resonance
in Nearly- Hamiltonian Systems, \textit{Stochastics and Dynamics},
\textbf{8}, 1, pp 1-21, 2008.

[4] Brin M., Freidlin M., On stochastic behavior of perturbed
Hamiltonian systems, \textit{Ergodic Theory and Dynamical Systems},
\textbf{20}, pp. 55 - 76, 2000.

[5] Bertotti G, Mayergoyz I., Serpico C., \textit{Nonlinear
Magnetization Dynamics in Nanosystems}, Elsevier, 2009.

[6] Dolgopyat D., Freidlin M., Koralov L., Deterministic and
Stochastic perturbations of area preserving flows on a
two-dimensional torus, \textit{Ergodic Theory and Dynamical
Systems}, to appear.

[7] Dolgopyat D., Koralov L., Averaging of incompressible flows on
two dimensional surfaces, preprint.

[8] Freidlin M., Sublimiting Distributions and Stabilization of
Solutions of Parabolic Equations with a Small Parameter,
\textit{Soviet Math. Dokl.}, \textbf{235}, 5, pp 1042-1045, 1977.

[9] Freidlin M., Metastability and stochastic resonance for
multiscale systems, \textit{Contemporary mathematics}, Volume
\textbf{469} (2008), pp 208-225.

[10] Freidlin M., Quasi-deterministic Approximation, Metastability
and Stochastic Resonance, \textit{Physica D}, \textbf{137}, pp
333-352, 2000.

[11] Freidlin M., Weber M., A remark on random perturbations of
nonlinear pendulum, \textit{Ann.Appl.Prob}, Vol. \textbf{9}, 1999,
No.3, 611-628.

[12] Freidlin M., Weber M., Random perturbations of dynamical
systems and diffusion processes with conservation laws,
\textit{Probab. Theory Relat. Fields}, \textbf{128}, 441-466 (2004).

[13] Freidlin M., Wentzell A., Random Perturbations of Hamiltonian
Systems, \textit{Mem. of AMS}, \textbf{523}, 1994.

[14] Freidlin M., Wentzell A., \textit{Random Perturbations of
Dynamical Systems}, Second edition, Springer, 1998.

[15] Freidlin M., Wentzell A., Diffusion processes on an open book
and the averaging principle, \textit{Stochastic Processes and their
Applications}, \textbf{113} (2004), 101-126.

[16] Freidlin M., Wentzell A., Long-time behavior of weakly coupled
oscillators, \textit{Journal of Statistical Physics}, Vol.
\textbf{123}, No. 6, June 2006.

[17] Hartman P., \textit{Ordinary Differential Equations}, John
Wiley and Sons, Inc., 1964.

[18] Landau L., Lifshitz E., On the theory of dispersion of magnetic
permeability in ferromagnetic bodies, in \textit{Collected Papers of
L.D.Landau}, pp. 101 - 114, Pergamon Press, 1965.

[19] Nikolaev I., Zhuzhoma E., \textit{Flows on 2-dimensional
manifolds. An overview}, Lect. Notes. Math., \textbf{1705} (1999),
Springer-Verlag, Berlin.

[20] Oksendal B., \textit{Stochastic Differential Equations}, Fifth
edition, Springer.

[21] Oliviery E., Vares M.E., \textit{Large Deviations and
Metastability}, Cambridge University Press, 2005.

\end{document}